\documentclass[reqno,oneside,12pt]{amsart}

\hoffset - 1.5cm

\usepackage{amsmath}
\usepackage{amsthm}
\usepackage{amssymb}
\usepackage{amstext}
\usepackage{amsopn}
\usepackage{amsfonts}

\newtheorem{Corollary}{Corollary}[section]

\newtheorem{Fact}{Fact}[section]
\newtheorem{Lemma}{Lemma}[section]
\newtheorem{Proposition}{Proposition}[section]
\newtheorem{Theorem}{Theorem}[section]
\newtheorem{Hipothesis}{Hipothesis}[section]

\theoremstyle{definition}
\newtheorem{Definition}{Definition}[section]
\newtheorem{Remark}{Remark}[section]
\newtheorem{Example}{Example}[section]

\newcommand{\ba}{\begin{array}}
\newcommand{\bc}{\begin{center}}
\newcommand{\bd}{\begin{description}}
\newcommand{\bdm}{\begin{displaymath}}
\newcommand{\be}{\begin{enumerate}}
\newcommand{\beq}{\begin{equation}}
\newcommand{\bdf}{\begin{Definition}}
\newcommand{\bex}{\begin{Example}}
\newcommand{\bft}{\begin{Fact}}
\newcommand{\bl}{\begin{Lemma}}
\newcommand{\bp}{\begin{Proposition}}
\newcommand{\br}{\begin{Remark}}
\newcommand{\bt}{\begin{Theorem}}
\newcommand{\bco}{\begin{Corollary}}
\newcommand{\bh}{\begin{Hipothesis}}
\newcommand{\ea}{\end{array}}
\newcommand{\ec}{\end{center}}
\newcommand{\ed}{\end{description}}
\newcommand{\edm}{\end{displaymath}}
\newcommand{\ee}{\end{enumerate}}
\newcommand{\eeq}{\end{equation}}
\newcommand{\edf}{\end{Definition}}
\newcommand{\eex}{\end{Example}}
\newcommand{\eft}{\end{Fact}}
\newcommand{\el}{\end{Lemma}}
\newcommand{\ep}{\end{Proposition}}
\newcommand{\er}{\end{Remark}}
\newcommand{\et}{\end{Theorem}}
\newcommand{\eco}{\end{Corollary}}
\newcommand{\eh}{\end{Hipothesis}}

\newcommand{\bH}{\mathbb{H}}
\newcommand{\bI}{\mathbb{I}}

\newcommand{\bK}{\mathbb{K}}

\newcommand{\bN}{\mathbb{N}}

\newcommand{\bR}{\mathbb{R}}

\newcommand{\bV}{\mathbb{V}}
\newcommand{\bW}{\mathbb{W}}

\newcommand{\bZ}{\mathbb{Z}}

\newcommand{\cC}{\mathcal{C}}

\newcommand{\cN}{\mathcal{N}}

\newcommand{\cU}{\mathcal{U}}

\newcommand{\numsec}{\setcounter{Theorem}{0}\setcounter{Definition}{0}
\setcounter{Remark}{0} \setcounter{Lemma}{0} \setcounter{Fact}{0}
\setcounter{Proposition}{0} \setcounter{Corollary}{0}
\setcounter{Example}{0} \setcounter{equation}{0}
\setcounter{Property}{0}\renewcommand\theequation{\arabic{section}.\arabic{equation}}
\renewcommand\theTheorem{\arabic{section}.\arabic{Theorem}}
\renewcommand\theDefinition{\arabic{section}.\arabic{Definition}}
\renewcommand\theRemark{\arabic{section}.\arabic{Remark}}
\renewcommand\theLemma{\arabic{section}.\arabic{Lemma}}
\renewcommand\theFact{\arabic{section}.\arabic{Fact}}
\renewcommand\theProposition{\arabic{section}.\arabic{Proposition}}
\renewcommand\theCorollary{\arabic{section}.\arabic{Corollary}}
\renewcommand\theExample{\arabic{section}.\arabic{Example}}
\renewcommand\theProperty{\arabic{section}.\arabic{Property}}}

\numberwithin{equation}{section} \errorcontextlines=0
\newcommand{\im}{\mathrm{ im \;}}
\newcommand{\diag}{\mathrm{ diag \;}}

\newcommand{\sign}{\mathrm{ sign \;}}
\newcommand{\dist}{\mathrm{ dist \;}}

\newcommand{\sone}{SO(2)}
\newcommand{\ds}{\displaystyle}
\newcommand{\nt}{\noindent}

\newcommand{\h}{\mathbb{H}}
\newcommand{\dg}{\nabla_{\sone}\mathrm{-deg}}
\newcommand{\ind}{\mathrm{ind}}
\newcommand{\dls}{\mathrm{deg}_{\mathrm{LS}}}

\begin{document}

\title[Bifurcation From Infinity]{Periodic Solutions  of  Second Order Hamiltonian \\
Systems  Bifurcating from Infinity \\ \vspace{0.3cm} Les solutions p\'eriodiques, \' emanants de l'infini des
syst\`emes hamiltoniens autonommes de second ordre}

\author{Justyna Fura$^{\dag}$}
\address{
Faculty of Mathematics and Computer Science\\
Nicolaus Copernicus University \\
PL-87-100 Toru\'{n} \\ ul. Chopina 12/18 \\
Poland} \email{Justyna.Fura@mat.uni.torun.pl}

\author{S{\l}awomir Rybicki$^{\ddag}$}
\address{Faculty of Mathematics and Computer Science\\
Nicolaus Copernicus University \\
PL-87-100 Toru\'{n} \\ ul. Chopina 12/18 \\
Poland} \email{Slawomir.Rybicki@mat.uni.torun.pl}

\date{\today}
\keywords{Autonomous second order Hamiltonian  systems; existence and continuation of periodic solutions; degree for
$\sone$-equivariant gradient maps;}

\subjclass[2000]{Primary: 37J45; Secondary: 37G15, 37G35.}

\thanks{$^{\dag}$ Partially sponsored by the Doctoral Program in Mathematics at the Nicolaus Copernicus University,
Toru\'n, Poland, and by the Scholarship of the Institute of Mathematics of the Polish Academy of Sciences}

\thanks{$^{\ddag}$Partially supported by the Ministry Education and Science; under grant 1 PO3A 009 27}

\begin{abstract}
The goal of this article is to study closed connected sets of periodic solutions,  of autonomous second order
Hamiltonian systems, emanating from infinity. The main idea is to apply the degree for $\sone$-equivariant gradient
operators defined by the second author in \cite{[RYB1]}. Using the results due to Rabier \cite{[RAB]}  we show that
we cannot apply the Leray-Schauder degree to prove the main results of this article. It is worth pointing out that
since we study connected sets of solutions, we also cannot use the Conley index technique and the Morse theory.
\end{abstract}
\maketitle

\bc R\'{e}sum\'{e} \ec

{\footnotesize Le but de cette article est l'\'etude des ensembles ferm\'es et connexes de solutions p\'eriodiques,
\' emanant de l'infini, des syst\`emes hamiltoniens autonomes de second ordre. L'id\'ee principale consiste a
\`{a}ppliquer le degr\'e aux op\'erateurs de gradient  $\sone $-\'equivariants d\'{e}finis par le second auteur dans
\cite{[RYB1]}. Moyennant un r\'esultat de Rabier \cite{[RAB]}, on d\'emontre que l'on ne peut pas appliquer le
degr\'e de Leray-Schauder pour obtenir le r\'esultat principal de ce travail. Il est important de souligner que, vu
que l'on \'etudie des ensembles connexes de solutions, ni la technique de l'indice de Conley, ni la th\'eorie de
Morse ne peuvent \^etre appliqu\'{e}es ici.}

\numsec
\section{Introduction}
Consider the following family of autonomous second order Hamiltonian  systems
\begin{equation}\label{introsys}
   \begin{cases}
    \ddot{u}(t) =  -  \nabla_u V(u(t),\lambda),&  \\
    u(0)=u(2\pi), & \\
    \dot{u}(0)=\dot{u}(2\pi),
  \end{cases}
\end{equation}
where $V \in C^2(\bR^n \times \bR,\bR)$ and the gradient $\nabla_x V$ (with respect to the first coordinate) is
asymptotically linear at infinity, i.e. $\nabla_x V(x,\lambda)=A(\lambda)x+o(\|x\|)$ as $\|x\| \rightarrow \infty$
uniformly on bounded $\lambda$-intervals and $A(\lambda)$ is a real symmetric matrix for every $\lambda\in\bR$.

Our purpose is to prove sufficient conditions for the existence of closed connected sets of non-stationary
$2\pi$-periodic solutions of system \eqref{introsys}    emanating from infinity.  Moreover, we describe the possible
minimal periods of solutions bifurcating from infinity and study the symmetry-breaking of solutions.

Bifurcations from infinity of     solutions of second order ODE's have been studied among the others in
\cite{[GLOVER1],[MA],[MALAGUTI],[SABATINI1],[SABATINI]}.  The authors applied the idea of the Hopf bifurcation from
infinity or the Leray-Schauder degree to study  solutions of the Li\'{e}nard, Rayleigh and Sturm-Liouville
equations. The assumptions considered in those articles are of different nature than these   in our article. For
example in the case of the Hopf bifurcation they considered asymptotically linear equation of the form $\ddot
x(t)=A(\lambda)x(t)+a(x,\lambda),$ where the matrix $A(\lambda)$ has a simple eigenvalue $i\omega_0$ ($0 \neq
\omega_0 \in \bR$) at $\lambda = \lambda_0$ and $a(x,\lambda) \rightarrow 0$ as $x \rightarrow 0$ i.e. the matrix
$A(\lambda_0)$ is not symmetric. Moreover, they do not obtain any estimation of minimal periods of bifurcating
solutions and information about the symmetry-breaking phenomenon.

We treat solutions of system \eqref{introsys} as critical orbits of an $\sone$-invariant $C^2$-functional $\Phi_V : \h^1_{2\pi} \times \bR
\rightarrow \bR$ whose gradient (with respect to the first coordinate) is an $\sone$-equivariant $C^1$-operator of the form compact perturbation
of the identity.

The basic idea is to apply the degree for $\sone$-equivariant gradient maps defined and discussed in
\cite{[RYB1],[RYB2],[RYB3],[RYB4]}. Our degree is an element of the tom Dieck ring $U(\sone),$ see Section
\ref{prelim} for the definition of this ring. The first degree for $\sone$-equivariant gradient maps, which is a
rational number, is due to Dancer \cite{[DANCER0]}. The degree for equivariant gradient maps in the presence of
symmetries of any compact Lie group $G$, which is an element of the tom Dieck ring $U(G)$, is due to G\c{e}ba
\cite{[GEBA0]}, see \cite{[DIECK]} for the definition of $U(G).$

For other applications of the degree for $\sone$-equivariant gradient maps to Hamiltonian  systems we refer the
reader to \cite{[FRR],[MARARY],[RADZKI],[RARY]}.

It is worth in pointing out that application of classical invariants like the Conley index technique and the Morse
theory does not ensure the existence of closed connected sets of critical points of variational problems, see
\cite{[AMB],[BOHME],[IZE0],[MARINO],[TAKENS]} for examples and discussion.

Since the gradient of the functional $\Phi_V$ is of the form compact perturbation of the identity, it is natural to
try to relate the degree for $\sone$-equivariant gradient maps to the Leray-Schauder degree. We are aware of
theorems similar to Theorem \ref{cont1} which have been proved for operators of the form compact perturbation of the
identity (without gradient and equivariant structures), see for instance Theorem 2.6 of \cite{[LESC]}.

However the choice of the degree for $\sone$-equivariant gradient maps seems to be the best adapted to our theory.
The advantage of using the degree for $\sone$-equivariant gradient maps lies in the fact that the index of an
isolated nontrivial $\sone$-orbit can be a nonzero element of the tom Dieck ring $U(\sone).$ Whereas the index of
this orbit computed by the Leray-Schauder degree equals $0 \in \bZ,$ see   \cite{[RAB]} .

After this introduction our article is organized as follows.

In Section \ref{prelim}, for the convenience of the reader, we have summarized without proofs the relevant material on the degree for
$\sone$-equivariant gradient maps, thus making our exposition self-contained.

Section \ref{continuation} is devoted to the study of closed connected sets of critical $\sone$-orbits of asymptotically linear
$\sone$-equivariant gradient maps of the form compact perturbation of the identity. Using the degree for $\sone$-equivariant gradient maps we
define a bifurcation index $\textsc{Bif}(\infty,[\lambda_-,\lambda_+])\in U(\sone)$, see Definition \ref{indinfty}. Nontriviality of the
bifurcation index implies the existence of an unbounded closed connected set of critical $\sone$-orbits, see Theorem \ref{cont1}. If the set of
stationary solutions of second order Hamiltonian  system is bounded then the bifurcation index
$\textsc{Bif}_{\textrm{LS}}(\infty,[\lambda_-,\lambda_+])\in \bZ$ computed by the Leray-Schauder degree is trivial. We discuss this situation in
Remarks \ref{rema}, \ref{uw3wn1cont1} and Corollary \ref{wn1cont1}. In Theorem \ref{cont2} we indicate points at which an unbounded closed
connected set of critical $\sone$-orbits meets infinity. In Lemmas \ref{lem1cont1}, \ref{lem1cont2} we control the isotropy groups of
$\sone$-orbits. The phenomenon of symmetry-breaking of $\sone$-orbits is discussed in Corollaries \ref{wn3cont2}, \ref{wn4cont2}.

In Section \ref{results} the main results of this article are stated and proved. In this section we study closed connected sets of periodic
solutions of autonomous second order Hamiltonian systems. Theorems \ref{eqcont1}, \ref{eqcont2} are consequences of Theorems \ref{cont1},
\ref{cont2}, respectively. In these theorems we have formulated sufficient conditions for the existence of  unbounded  closed connected sets of
$2\pi$-periodic solutions of system \eqref{introsys}.  We emphasize that assumptions of these theorems are expressed directly in terms of the
right hand sight of system \eqref{introsys} i.e. potential $V.$ In Corollary \ref{wn1eqcont2} we have described the minimal periods of solutions
of system \eqref{introsys} which are sufficiently close to infinity. In Theorem \ref{eqcont3} we study periodic solutions of a special case of
system \eqref{introsys} i.e. we assume that $V(x,\lambda)=\lambda^2 V(x).$ In this theorem we indicate all the points at which closed connected
sets of periodic solutions of system \eqref{introsys} meet infinity. The minimal periods of solutions of system  \eqref{introsys} are discussed in
Corollary \ref{wn1eqcont3}.

In Section \ref{examples} we consider three real second order Hamiltonian  systems in order to illustrate the main
results of this paper.

\nt   \textsc{Acknowledgments.} The authors wishes to express their thanks to the referee  for several helpful
comments concerning the style and exposition of this article.

\numsec
\section{Preliminaries}
\label{prelim}

In this section, for the convenience of the reader, we remind the main properties of the degree for
$\sone$-equivariant gradient maps defined in \cite{[RYB1]}. This degree will be denoted briefly by
$\nabla_{\sone}\mathrm{-deg}$ to underline that it is a special degree theory for $\sone$-equivariant gradient maps.

Put $\displaystyle U(\sone)= \bZ \oplus   \bigoplus_{k=1}^{\infty} \bZ$ and define the actions
$$+ , \star : U(\sone) \times U(\sone) \rightarrow U(\sone),$$
$$\cdot : \bZ \times U(\sone) \rightarrow U(\sone),$$ as follows

\begin{align}
  \label{doda} \alpha + \beta =&\left(\alpha_0 + \beta_0, \alpha_1 + \beta_1,
 \ldots,\alpha_k+\beta_k,\ldots\right),   \\
 \label{m} \alpha \star \beta =&(\alpha_0   \beta_0, \alpha_0  \beta_1 +
\beta_0  \alpha_1,   \ldots, \alpha_0  \beta_k + \beta_0  \alpha_k, \ldots), \\
\gamma \cdot \alpha = &(\gamma  \alpha_0, \gamma  \alpha_1,   \ldots, \gamma  \alpha_k, \ldots),
\end{align}
where $\alpha = (\alpha_0, \alpha_1, \ldots, \alpha_k, \ldots), \beta = (\beta_0, \beta_1, \ldots,\beta_k, \ldots)
\in U(\sone)$ and $\gamma \in \bZ.$ It is easy to check that $(U(\sone),+,\star)$ is a commutative ring with  the
trivial element $\Theta=(0,0,\ldots) \in U(\sone)$ and the unit $\bI=(1,0,\ldots) \in U(\sone).$  The ring
$(U(\sone),+,\star)$ is called the tom Dieck ring of the group $\sone$. For the definition of the tom Dieck ring
$U(G),$ where $G$ is any compact Lie group, we refer the reader to \cite{[DIECK]}.

\nt If $\delta_1,\ldots,\delta_q \in U(\sone),$ then we write $\ds \prod_{j=1}^q \delta_j$ for $\delta_1 \star
\ldots \star \delta_q.$ Moreover, it is understood that $\ds \prod_{j \in \emptyset} \delta_j= \bI \in U(\sone).$

A representation of the group $\sone$  (an $\sone$-representation) is a pair $\bV=(\bV_0,\rho),$ where $\bV_0$ is a
real, linear space and  $\rho : \sone \rightarrow GL(\bV_0)$ is a continuous homomorphism into the group of all
linear automorphisms of $\bV_0.$ Notice that if $\bV=(\bV_0,\rho)$ is an $\sone$-representation, then letting
$gv=\rho(g)(v)$ we obtain a linear $\sone$-action on $\bV_0.$ For simplicity of notation, we do not distinguish
between $\bV$ and $\bV_0$ using the same letter $\bV$ for a representation and the corresponding linear space
$\bV_0$ .

Let $\bV$ be a real, finite-dimensional and orthogonal $\sone$-representation. If $v \in \bV$ then  the subgroup
$SO(2)_v=\{g \in \sone : g v =v\}$ is said to be the isotropy group of $v \in \bV.$ Moreover, the set $\sone v =\{gv
: g \in \sone\}$ is called the $\sone$-orbit of $v \in \bV.$ \\ Let $\Omega \subset \bV$ be an open, bounded and an
$\sone$-invariant subset and let $H \subset \sone$ be a closed subgroup. Then we define
\begin{itemize}
  \item $\Omega^H=\{v \in \Omega : H \subset \sone_v\}=\{v \in \Omega : g v = v \: \forall \: g \in H\},$
  \item $\Omega_H=\{v \in \Omega : H = \sone_v\}.$
\end{itemize}

Fix $k \in \bN$ and set $C^k_{\sone}(\bV,\bR) = \{f \in C^k(\bV,\bR) : f \text{ is } \sone\text{-invariant}\}.$

Let $f_0 \in C^1_{\sone}(\bV,\bR).$ Since $\bV$ is an orthogonal $\sone$-representation, the gradient  $\nabla f_0 :
\bV \rightarrow \bV$ is an $\sone$-e\-qui\-va\-riant $C^0$-map. If $H \subset \sone$ is a closed subgroup then
$\bV^H$ is a finite-dimensional $\sone$-representation. If $\big(\nabla f_0\big)^H=\nabla f_{0 \mid \bV^H}$ and
$f_0^H=f_{0 \mid \bV^H }$ then it is easy to verify that  $\big(\nabla f_0\big)^H =\nabla \big(f_0^H \big) : \bV^H
\rightarrow \bV^H$ is well-defined $\sone$-equivariant gradient map.

Choose an open, bounded and $\sone$-invariant subset $\Omega \subset \bV$ such that $(\nabla f_0)^{-1}(0) \cap
\partial \Omega = \emptyset.$ Under these assumptions we have defined in \cite{[RYB1]} the degree for
$\sone$-equivariant gradient maps $\nabla_{\sone}\mathrm{-deg}(\nabla f_0,\Omega) \in U(\sone)$ with coordinates
$$\nabla_{\sone}\mathrm{-deg}(\nabla f_0,\Omega)=$$$$=(\nabla_{\sone}\mathrm{-deg}_{\sone}(\nabla f_0,\Omega),
\nabla_{\sone}\mathrm{-deg}_{\bZ_1}(\nabla f_0,\Omega), \ldots, \nabla_{\sone}\mathrm{-deg}_{\bZ_k}(\nabla
f_0,\Omega), \ldots ).$$

\br \label{goodrem} To define the degree for $\sone$-equivariant gradient maps of $\nabla f_0$ we choose (in a
homotopy class of the $\sone$-equivariant gradient map $\nabla f_0$) a "sufficiently good" $\sone$-equivariant
gradient map $\nabla f_1$ and define  this degree for $\nabla f_1.$ The definition does not depend on the choice of
the map $\nabla f_1.$ Roughly speaking the main steps of the definition of the degree for $\sone$-equi\-va\-riant
gradient maps of $\nabla f_0 : (cl(\Omega),\partial \Omega) \rightarrow (\bV, \bV \setminus \{0\})$ are the
following:
\begin{enumerate}
\item[Step 1.] There is a potential  $f \in C^1_{\sone}(\bV \times [0,1],\bR)$ such that
\begin{enumerate}
\item[(a1)] $(\nabla_v f)^{-1}(0) \cap (\partial \Omega \times [0,1])=\emptyset,$
\item[(a2)] $\nabla_v f(\cdot,0)=\nabla f_0(\cdot),$
\item[(a3)] $\nabla_v f_1 \in C^1_{\sone}(\bV,\bV),$ where we abbreviate $\nabla_v f(\cdot,1)$ to $\nabla_v f_1,$
\item[(a4)] $(\nabla_v f_1)^{-1}(0) \cap \Omega^{\sone}=\{v_1,\ldots,v_p\}$ and
\begin{enumerate}
\item[(i)] $\det \nabla^2_{vv} f_1(v_j) \neq 0,$ for all $j =1,\ldots,p,$
\item[(ii)] $\nabla^2_{vv} f_1(v_j)=\left[\begin{array}{cc}
  \nabla^2_{vv}\big(f_1^{\sone}\big)(v_j) & 0 \\
  0 & Id
\end{array} \right] :   \begin{array}{c}
  \bV^{\sone} \\ \oplus \\
  (\bV^{\sone})^{\perp}
\end{array} \longrightarrow  \begin{array}{c}
  \bV^{\sone} \\ \oplus \\
  (\bV^{\sone})^{\perp}
\end{array},$ for all $j =1,\ldots,p,$
\end{enumerate}
\item[(a5)] $(\nabla_v f_1)^{-1}(0) \cap (\Omega \setminus \Omega^{\sone})=\{\sone w_1, \ldots, \sone w_q\}$
and
\begin{enumerate}
\item[(i)] $\dim \ker \nabla^2_{vv} f_1(w_j)=1,$ for all $j =1,\ldots,q,$
\item[(ii)] $$\nabla^2_{vv} f_1(w_j)= \left[\begin{array}{ccc}
  0 & 0 &  0\\
  0 & Q_j & 0 \\
  0 & 0 & Id
\end{array}\right]   :$$ $$ \begin{array}{c}
  T_{w_j} (SO(2)w_j) \\ \oplus \\ T_{w_j}(\bV_{\sone_{w_j}}) \ominus  T_{w_j} (SO(2)w_j)\\ \oplus \\
 ( T_{w_j}(\bV_{\sone_{w_j}}))^{\perp}
\end{array}   \longrightarrow  \begin{array}{c}
  T_{w_j} (SO(2)w_j) \\ \oplus \\ T_{w_j}(\bV_{\sone_{w_j}}) \ominus  T_{w_j} (SO(2)w_j)\\ \oplus \\
 ( T_{w_j}(\bV_{\sone_{w_j}}))^{\perp}
\end{array},$$ for all $j =1,\ldots,q.$
\end{enumerate}
\end{enumerate}

\item[Step 2.] The first coordinate of the degree for $\sone$-equivariant gradient maps is defined by
$\ds \dg_{\sone}(\nabla f_0,\Omega)=\sum_{j=1}^p \sign \det \nabla^2_{vv}(f_1^{\sone})(v_j).$ In other words since
$\nabla \big(f_1^{\sone} \big) = \big(\nabla f_1\big)^{\sone},$ we obtain $$\dg_{\sone}(\nabla f_0,\Omega)=
\mathrm{deg}_{\mathrm{B}}((\nabla f_1)^{\sone},\Omega^{\sone},0),$$ where $\mathrm{deg}_{\mathrm{B}}$ denotes  the
Brouwer degree.
\item[Step 3.] Fix $k \in \bN$ and define $$\dg_{\bZ_k}(\nabla f_0,\Omega)=
\sum_{\{j\in\{1,\ldots,q\} : \sone_{w_j}=\bZ_k\}} \sign \det Q_j,$$
\end{enumerate}
Notice that since $$\mathrm{deg}_{\mathrm{B}}((\nabla
f_1)^{\sone},\Omega^{\sone},0)=\mathrm{deg}_{\mathrm{B}}(\nabla f_1,\Omega,0) \textrm{ and }
\mathrm{deg}_{\mathrm{B}}(\nabla f_1,\Omega,0)=\mathrm{deg}_{\mathrm{B}}(\nabla f_0,\Omega,0)$$ (see \cite{[RAB]}),
directly by the Step 2.   we obtain $\dg_{\sone}(\nabla f_0,\Omega)= \mathrm{deg}_{\mathrm{B}}(\nabla
f_0,\Omega,0).$ Moreover, immediately from the Step 3. we obtain that if $k \in \bN$ and $\sone_v \neq \bZ_k$ for
every $v \in \Omega,$ then $\dg_{\bZ_k}(\nabla f_0,\Omega)=0.$ \er

For $\gamma > 0$ and $v_0 \in \bV^{\sone}$ we put $B_{\gamma}(\bV,v_0) = \{v \in \bV :\ \mid v - v_0 \mid <
\gamma\}.$ For simplicity of notation, we write $B_{\gamma}(\bV)$ instead of $B_{\gamma}(\bV,0).$ In the following
theorem we formulate the main properties of the degree for $\sone$-equivariant gradient maps.

\bt[\cite{[RYB1]}]\label{wlas} Under the above assumptions the degree for $\sone$-e\-qui\-va\-riant gradient maps
has the following properties\et

\begin{enumerate}
  \item  \label{w1}   if
$\nabla_{\sone}\mathrm{-deg}(\nabla f, \Omega) \neq \Theta,$ then $(\nabla f)^{-1}(0) \cap \Omega \neq \emptyset,$
\item  if  $\nabla_{\sone}\mathrm{-deg}_H(\nabla f, \Omega) \neq 0,$ then $(\nabla f)^{-1}(0) \cap \Omega^H
\neq \emptyset,$
  \item  \label{w3} if
$\Omega = \Omega_0 \cup \Omega_1$ and $\Omega_0 \cap \Omega_1 = \emptyset,$ then
$$\nabla_{\sone}\mathrm{-deg}(\nabla f, \Omega) = \nabla_{\sone}\mathrm{-deg}(\nabla f, \Omega_0) +
\nabla_{\sone}\mathrm{-deg}(\nabla f, \Omega_1),$$
  \item  \label{w4}  if
$ \Omega_0 \subset \Omega$ is an open $\sone$-invariant subset and $(\nabla f)^{-1}(0) \cap \Omega \subset
\Omega_0,$ then
$$ \nabla_{\sone}\mathrm{-deg}(\nabla f, \Omega) = \nabla_{\sone}\mathrm{-deg}(\nabla f, \Omega_0),$$
\item   \label{w2} if    $ f \in C^1_{\sone}(\bV \times [0,1],\bR)$  is such that
$(\nabla_v f)^{-1}(0) \cap \left(\partial \Omega \times [0,1] \right) = \emptyset,$ then
$$\nabla_{\sone}\mathrm{-deg}(\nabla f_0, \Omega)=\nabla_{\sone}\mathrm{-deg}(\nabla f_1,\Omega),$$

\item \label{wzaw}  if $\bW$ is an orthogonal $\sone$-representation,     then
  $$\nabla_{\sone}\mathrm{-deg}((\nabla f, Id), \Omega \times B_{\gamma}(\bW)) =
\nabla_{\sone}\mathrm{-deg}(\nabla f, \Omega),$$
\item  if  $ f \in C^2_{\sone}(\bV,\bR)$ is such that $\nabla f(0) = 0$ and $\nabla^2f(0)$
is an   $\sone$-equivariant self-adjoint isomorphism  then there is $\gamma > 0$ such that
$$\nabla_{\sone}\mathrm{-deg}(\nabla f, B_{\gamma}(\bV)) = \nabla_{\sone}\mathrm{-deg}(\nabla^2 f(0), B_{\gamma}(\bV)).$$
\end{enumerate}

\nt Below we formulate the product formula for the degree for $\sone$-equivariant gradient maps.

\bt[\cite{[RYB2]}] \label{pft} Let $\Omega_i \subset \bV_i$ be an open, bounded and $\sone$-invariant subset of a
finite-dimensional, orthogonal $\sone$-representation $\bV_i,$ for $i=1,2.$ Let $f_i \in C^1_{\sone}(\bV_i,\bR) $ be
such that $\big(\nabla f_i \big)^{-1}(0) \cap
\partial \Omega_i = \emptyset,$  for $i=1,2.$  Then
$$
\nabla_{\sone}\mathrm{-deg}((\nabla f_1, \nabla f_2), \Omega_1
\times \Omega_2) = \nabla_{\sone}\mathrm{-deg}(\nabla f_1,
\Omega_1) \star \nabla_{\sone}\mathrm{-deg}(\nabla f_2, \Omega_2).
$$
\et

\noindent For $k \in \bN$ define a map $\rho^k : \sone \rightarrow GL(2,\bR)$ as follows
\[
\rho^k (e^{i\theta})= \left[\begin{array}{lr}
\cos (k\theta)&-\sin (k\theta)\\
\sin (k\theta)&\cos (k\theta)
\end{array}\right]
\qquad 0\le\theta  <  2  \pi.
\]
For $j,k \in \bN$ we denote by $\bR[j,k]$ the direct sum of $j$ copies of $(\bR^2 ,\rho^k)$, we also denote by
$\bR[j,0]$ the trivial $j$-dimensional $\sone$-representation. We say that two $\sone$-representations $\bV$ and
$\bW$ are equivalent  if there exists an $\sone$-equivariant, linear isomorphism $T : \bV \rightarrow \bW$. The
following classic result gives  complete classification (up to equivalence) of finite-di\-men\-sio\-nal
representations of the group $\sone$ (see \cite{[ADM]}).

\begin{Theorem}[\cite{[ADM]}]
\label{tk} If $\bV$ is a finite-dimensional $\sone$-representation,  then there exist finite sequences $\{j_i\},\,
\{k_i\}$ satisfying:\\ $ (*)\qquad k_i\in \{0\}\cup \bN,\quad  j_i\in \bN,\quad  1\le i\le r,  \: k_1 < k_2 < \dots
< k_r $
\\  such that $\bV$ is equivalent to $\displaystyle\bigoplus^r_{i=1} \bR[j_i ,k_i]$. Moreover, the
equivalence class of $\bV$, ($\bV\approx\displaystyle\bigoplus^r_{i=1} \bR[j_i,k_i]$) is uniquely determined by
$\{k_i\},\, \{j_i\}$ satisfying $(*)$.
\end{Theorem}

\nt We will denote by $m^-(L)$ the Morse index of a symmetric matrix $L$ i.e. the sum of algebraic multiplicities of negative eigenvalues of $L.$

\nt To apply successfully any degree theory we need computational formulas for this invariant. Below we show how to
compute the degree for $\sone$-equivariant gradient maps of a linear, self-adjoint, $\sone$-equivariant isomorphism.

\bl[\cite{[RYB1]}] \label{lindeg} If $\bV \approx \bR[j_0,0] \oplus \bR[j_1,k_1] \oplus \ldots \oplus \bR[j_r,k_r],$
 $L : \bV \rightarrow \bV$ is a self-adjoint, $\sone$-equivariant, linear isomorphism and $\gamma > 0$ then \el
\begin{enumerate}
  \item $L= \diag (L_0,L_1,\ldots,L_r),$
  \item $$\nabla_{\sone}\mathrm{-deg}_H(L,B_{\gamma}(\bV))=
 \begin{cases}
    (-1)^{m^-(L_0)}, & \text{ for } H = \sone, \\
   \displaystyle (-1)^{m^-(L_0)}  \cdot
     \frac{m^-(L_i)}{2}, & \text{ for } H = \bZ_{k_i}\\
     0, & \text{ for } H \notin \{\sone, \bZ_{k_1}, \ldots, \bZ_{k_r}\},
  \end{cases}
$$
  \item in particular, if $L=-Id,$ then
 $$\nabla_{\sone}\mathrm{-deg}_H(-Id,B_{\gamma}(\bV))=
 \begin{cases}
    (-1)^{j_0}, & \text{ for } H = \sone, \\
   \displaystyle (-1)^{j_0}  \cdot j_i, & \text{ for } H = \bZ_{k_i},\\
     0, & \text{ for } H \notin \{\sone, \bZ_{k_1}, \ldots, \bZ_{k_r}\}.
  \end{cases}
$$
\end{enumerate}

\noindent Let $(\bH, \langle \cdot,\cdot \rangle _{\bH})$ be an infinite-dimensional, separable Hilbert space which
is an or\-tho\-go\-nal $\sone$-representation  and let $C_{\sone}^k(\bH ,\bR)$ denote the set of $\sone$-invariant
$C^k$-functionals. Fix $\Phi \in C_{\sone}^1(\bH ,\bR)$ such that $\nabla \Phi(u)=   u -  \nabla \eta(u),$ where
$\nabla  \eta : \bH   \rightarrow \bH$ is an $\sone$-equivariant   compact operator. Let $\cU \subset \h$ be an
open, bounded and $\sone$-invariant set such that $\left(\nabla \Phi \right)^{-1}(0) \cap \partial \cU = \emptyset.$
In this situation $\ds \dg(Id - \nabla \eta, \cU) \in U(\sone)$ is well-defined, see \cite{[RYB1]} for details  and
basic properties of this degree.

\br \label{nowyrem}  We would like to underline that the infinite-dimensional version of the degree for
$\sone$-equivariant gradient maps has the following two important properties
\begin{enumerate}
  \item $\ds \dg_{\sone}(\nabla \Phi, \cU)=\mathrm{deg}_{\mathrm{LS}}(\nabla \Phi, \cU,0),$ where
$\mathrm{deg}_{\mathrm{LS}}$ denotes the Leray-Schau\-der degree,
  \item if $\Phi \in C^1_{\sone}(\bH \times [\lambda_-,\lambda_+],\bR), Q \subset \bH \times [\lambda_-,\lambda_+]$
  is an open bounded $\sone$-in\-va\-riant subset and   there is $\gamma > 0$ such that
\begin{enumerate}
  \item[a)]
  $Q \cap (\bH \times \{\lambda_-,\lambda_+\})=B_{\gamma}(\bH) \times \{\lambda_-,\lambda_+\},$
  \item[b)] $(\nabla_u \Phi)^{-1}(0) \cap \partial Q \subset B_{\gamma}(\bH) \times \{\lambda_-,\lambda_+\},$
\end{enumerate}
then $\dg(\nabla_u \Phi(\cdot,\lambda_+),Q_{\lambda_+})=\dg(\nabla_u \Phi(\cdot,\lambda_-),Q_{\lambda_-}),$ where
$Q_{\lambda_{\pm}}=\{(u,\lambda_{\pm}) \in Q\}.$
\end{enumerate}
The second property is a slight generalization of the homotopy invariance of the degree for $\sone$-equivariant
gradient maps and is called the generalized homotopy invariance. \er

\nt Let $L : \h \rightarrow \h$ be a linear, bounded, self-adjoint, $\sone$-equivariant operator with spec\-trum
$\sigma(L)=\{\lambda_i\}.$ By $\bV_L(\lambda_i)$ we will denote the eigenspace of $L$ corresponding to the
e\-igen\-va\-lue $\lambda_i$ and we put $\mu_L (\lambda_i)=\dim \bV_L(\lambda_i).$ In other words $\mu_L(\lambda_i)$
is the multiplicity of the eigenvalue $\lambda_i.$ Since operator $L$ is linear, bounded, self-adjoint, and
$\sone$-equivariant, $\bV_L(\lambda_i)$ is a finite-dimensional, orthogonal $\sone$-representation. For $\gamma > 0$
and $v_0 \in \h^{\sone}$ set $B_{\gamma}(\h,u_0)=\{u\in\h : \| u - u_0 \| < \gamma\}$.  For abbreviation, let
$B_{\gamma}(\h)$ stand for $B_{\gamma}(\h,0).$ Note that $B_{\gamma}(\h,u_0)$ is  open and $\sone$-invariant  for
every $u_0 \in \h^{\sone}$.

\nt Combining Theorem 4.5 in  \cite{[RYB1]} with Theorem \ref{pft} we obtain the following theorem.

\bt  \label{dizom}  Under the above assumptions if $1 \notin \sigma(L),$ then $$\dg(Id - L, B_{\gamma}(\h))=
\prod_{\lambda_i > 1} \dg(-Id, B_{\gamma}(\bV_L(\lambda_i))) \in U(\sone).$$ It is understood that if $\sigma(L)
\cap [1,+\infty)=\emptyset,$ then $$\dg(Id - L, B_{\gamma}(\h))=\bI \in U(\sone).$$ \et

\section{Abstract results}
\label{continuation}

\numsec

In this section we study   global bifurcation from infinity of critical orbits of $\sone$-invariant functionals.

\nt Let $(\h,\langle \cdot,\cdot \rangle_{\h})$ be as in the previous section. We consider $\h\times\bR$ as an
$\sone$-representation with $\sone$-action given by $g(u,\lambda)=(gu,\lambda),$ where $(u,\lambda) \in \h \times
\bR$ and $g \in \sone.$   Put $C^k_{\sone}(\h\times\bR,\bR)=\{\Phi\in C^k(\h\times\bR,\bR) : \Phi \;\text{is }
\sone\text{-invariant}\}.$ It is clear that if $\Phi\in C^k_{\sone}(\h\times\bR,\bR)$, then the gradient $\nabla_u
\Phi: \h\times\bR\rightarrow\h$ is an $\sone$-equivariant $C^{k-1}$-operator.

Consider a potential $\Phi \in C^2_{\sone}(\h\times\bR,\bR)$ such that:

\begin{enumerate}
\item[\textrm{\bf (c1)}] $\ds{\Phi(u,\lambda)=\frac 12 \langle u,u \rangle_{\h} - g(u,\lambda),}$
where $\nabla_u g:\h\times\bR\rightarrow\h$ is compact.
\end{enumerate}

\nt From now on we study solutions of the following system
\begin{equation}\label{rbif}
  \nabla_u \Phi(u,\lambda)  =0.
\end{equation}

\nt The set $(\nabla_u \Phi)^{-1}(0) \cap (\h^{\sone}\times\bR)$ is called the set of trivial solutions of equation
\eqref{rbif}. Put
$$\ds{\cN(\nabla_u \Phi)= \{(u,\lambda)\in (\h \setminus \h^{\sone}) \times\bR : \nabla_u
\Phi(u,\lambda)=0}\}.$$

Assume that there exist $\lambda_-, \lambda_+ > 0$ and $\gamma>0$ such that

\begin{equation}\label{zera}
(\nabla_u \Phi(\cdot,\lambda_{\pm}))^{-1}(0)\cap \left((\h\setminus
B_{\gamma}(\h))\times\{\lambda_{\pm}\}\right)=\emptyset.
\end{equation}

\bdf\label{indinfty} An element $\textsc{Bif}(\infty,[\lambda_-,\lambda_+])\in U(\sone)$ defined as follows
$$\textsc{Bif}(\infty,[\lambda_-,\lambda_+])=\dg(\nabla_u
\Phi(\cdot,\lambda_+),B_{\gamma}(\h))-\dg(\nabla_u \Phi(\cdot,\lambda_-),B_{\gamma}(\h))$$ is called the bifurcation
index at $(\infty,[\lambda_-,\lambda_+]).$ \edf

\nt The following lemma will be extremely useful in the proof  of the next theorem.

\bl(\cite{[Brown]})\label{sep} Let A and B be disjoint closed subsets of a compact space K. If there is no closed connected subset of K that
intersects both A and B, then there exist disjoint closed subsets $K_A$ and $K_B$ of K such that $A \subset K_A, B \subset K_B$ and $K=K_A \cup
K_B$. \el

\nt The following theorem is the most general result of this section. Namely, we prove the sufficient condition for
the existence of an unbounded closed connected set of critical orbits of $\sone$-invariant functionals. In the proof
of this theorem we combine Lemma \ref{sep} with the  degree for $\sone$-equivariant gradient maps.

\bt\label{cont1} Let $\Phi \in C^2_{\sone}(\h\times\bR,\bR)$ satisfy condition \textrm{\bf(c1)} and let
$\lambda_{\pm}\in\bR,\; \gamma>0$ be such that \eqref{zera} holds. If
$\textsc{Bif}(\infty,[\lambda_-,\lambda_+])\neq\Theta\in U(\sone),$ then there exists an unbounded closed connected
component $C$ of $(\nabla_u \Phi)^{-1}(0) \cap (\h \times [\lambda_-,\lambda_+])$ such that \linebreak $C \cap
\left(B_{\gamma}(\h)\times\{\lambda_-,\lambda_+\}\right) \neq \emptyset.$\et
\begin{proof}
First of all we claim that for every $\xi\geq \gamma$ there exists a closed connected component $C_{\xi} $ of
$\nabla_u \Phi^{-1}(0) \cap (\h \times [\lambda_-,\lambda_+])$ such that $$C_{\xi} \cap
\left(B_{\gamma}(\h)\times\{\lambda_-,\lambda_+\}\right)\neq \emptyset \textrm{  and } C_{\xi}\cap (\partial
B_{\xi}(\h)\times[\lambda_-,\lambda_+]) \neq \emptyset.$$ Suppose, contrary to our claim, that there exists $\xi
\geq \gamma$ such that at least one of the following conditions is fulfilled
\begin{enumerate}
  \item[(i)] $\cC \cap \left(B_{\gamma}(\h)\times\{\lambda_-,\lambda_+\}\right)=\emptyset,$
  \item[(ii)] $\cC \cap (\partial B_{\xi}(\h)\times[\lambda_-,\lambda_+])= \emptyset,$
\end{enumerate}
for every closed connected component $\cC$ of $\nabla_u \Phi^{-1}(0) \cap (\h \times [\lambda_-,\lambda_+])$.

\nt Put in Lemma \ref{sep}
\begin{enumerate}
\item[(i)] $\ds K= \nabla_u \Phi^{-1}(0) \cap
\left(cl(B_{\xi}(\h))\times[\lambda_-,\lambda_+]\right),$
\item[(ii)] $\ds A =\nabla_u \Phi^{-1}(0) \cap
\left(cl(B_{\xi}(\h))\times\{\lambda_-,\lambda_+\}\right),$
\item[(iii)] $\ds B= \nabla_u \Phi^{-1}(0) \cap \left(\partial
B_{\xi}(\h)\times[\lambda_-,\lambda_+]\right).$
\end{enumerate}

\nt Since $\nabla_u \Phi$ is of the form compact perturbation of the identity and
$cl(B_{\xi}(\h))\times[\lambda_-,\lambda_+]$ is closed and bounded, $K$ is compact. Recall that $\nabla_u
\Phi^{-1}(0) \cap \left(cl(B_{\xi}(\h))\times\{\lambda_-,\lambda_+\}\right)\subset
B_{\gamma}(\h)\times\{\lambda_-,\lambda_+\}.$ Thus $A \cap B=\emptyset$. By assumption, there is no closed connected
subset of $K$ that intersects both $A$ and $B$. Applying Lemma \ref{sep}, we obtain compact sets $K_A, K_B$ with
desired properties.

\nt Choose $\alpha>0$ such that $K_A(\alpha), K_B(\alpha)$ are disjoint $\alpha$-neighborhoods of the sets $K_A,
K_B$. Define $$\ds Q=\sone\left(\left(B_{\xi}(\h)\times[\lambda_-,\lambda_+]\right) \setminus
cl(K_B(\alpha))\right)=$$ $$ =\{(gv,\lambda) : v \in (B_{\xi}(\h)\times[\lambda_-,\lambda_+])  \setminus
cl(K_B(\alpha)) \textrm{ and } g \in \sone \}.$$

\nt We claim that $Q$ is open, $\sone$-invariant  and $(\nabla_u \Phi)^{-1}(0) \cap \partial Q \subset
B_{\xi}(\h)\times\{\lambda_-,\lambda_+\}.$ Since $B_{\xi}(\h)\times[\lambda_-,\lambda_+]$ is open in
$\h\times[\lambda_-,\lambda_+],$ it is clear that $Q$ is open. Moreover, since $Q$ is a sum of $\sone$-orbits, it is
$\sone$-invariant. What is left is to show that $(\nabla_u \Phi)^{-1}(0) \cap \partial Q \subset
B_{\xi}(\h)\times\{\lambda_-,\lambda_+\}.$ Suppose, contrary to our claim that, $(\nabla_u \Phi)^{-1}(0) \cap
(\partial Q \setminus (B_{\xi}(\h)\times\{\lambda_-,\lambda_+\})) \neq \emptyset$ and fix $(u_0,\lambda_0) \in
\partial Q \setminus (B_{\xi}(\h)\times\{\lambda_-,\lambda_+\})$ such that $\nabla_u \Phi(u_0,\lambda_0)=0.$ Hence
there are $(\widetilde{u}_0,\lambda_0) \in \partial ((B_{\xi}(\h)\times[\lambda_-,\lambda_+])  \setminus
cl(K_B(\alpha)))$ and $g \in \sone$ such that $(g\widetilde{u}_0,\lambda_0)=(u_0,\lambda_0).$ Since $\nabla \Phi$ is
$\sone$-equivariant, we obtain $$0=\nabla_u \Phi(u_0,\lambda_0)=\nabla_u \Phi(g \widetilde{u}_0,\lambda_0)=g
\nabla_u \Phi(\widetilde{u}_0,\lambda_0)$$ and consequently $\nabla_u \Phi(\widetilde{u}_0,\lambda_0)=0,$ which
contradicts the definition of $K_B(\alpha).$ \\ Put $\ds Q_{\lambda}=\{(u,\lambda) \in Q\}$ for every $\lambda \in
[\lambda_-,\lambda_+].$

Since  $(\nabla_u \Phi)^{-1}(0) \cap  \partial Q \subset B_{\gamma}(\h)\times\{\lambda_-,\lambda_+\}$, from the
generalized homotopy invariance of the degree for $\sone$-equivariant gradient maps (see Remark \ref{nowyrem}), we
obtain that:
\begin{align*}
\Theta &=\dg(\nabla_u \Phi(\cdot,\lambda_+),Q_{\lambda_+}) - \dg(\nabla_u \Phi(\cdot,\lambda_-),Q_{\lambda_-}) =
\\
&=  \dg(\nabla_u \Phi(\cdot,\lambda_+),B_{\xi}(\h)) - \dg(\nabla_u \Phi(\cdot,\lambda_-),B_{\xi}(\h)) =
\\
&=  \dg(\nabla_u \Phi(\cdot,\lambda_+),B_{\gamma}(\h)) - \dg(\nabla_u \Phi(\cdot,\lambda_-),B_{\gamma}(\h)) =
 \\ &= \textsc{Bif}(\infty,[\lambda_-,\lambda_+])\neq \Theta,
\end{align*}
\nt a contradiction.

Suppose, contrary to our claim that,  the theorem is false i.e. every closed connected component $C$ of $(\nabla_u
\Phi)^{-1}(0) \cap (\h \times [\lambda_-,\lambda_+])$ such that  $C \cap \left(B_{\gamma}(\h)\times
\{\lambda_-,\lambda_+\}\right)\neq\emptyset$ is bounded. Choose   an increasing sequence $\{\gamma_n\}\subset \bN$
such that $\gamma_n \geq \gamma$ for every $n \in \bN.$  From the first part of the proof it is known that  for
every $n\in\bN$ there exists a bounded closed connected component $C_{\gamma_n}$ of $(\nabla_u \Phi)^{-1}(0) \cap
(\h \times [\lambda_-,\lambda_+])$ such that $C_{\gamma_n} \cap
\left(B_{\gamma}(\h)\times\{\lambda_-,\lambda_+\}\right)\neq\emptyset$ and $C_{\gamma_n}\cap (\partial
B_{\gamma_n}(\h)\times[\lambda_-,\lambda_+])\neq \emptyset$. Choose $(u_n,\lambda_n)\in C_{\gamma_n} \cap
\left(B_{\gamma}(\h)\times\{\lambda_-,\lambda_+\}\right)$ for every $n\in \bN$. Without loosing of generality, one
can assume that $\lambda_n= \lambda_+$ for every $n\in \bN$. Note that $cl \{(u_n,\lambda_+)\}$ is compact, as a
closed subset of the compact set $(\nabla_u \Phi)^{-1}(0)\cap\left(cl(B_{\gamma}(\h))\times\{\lambda_+\}\right)$.
Thus, there exists convergent subsequence $(u_{n_k},\lambda_+) \rightarrow (u_0,\lambda_+).$ Denote by $C$ a closed
connected component of $(\nabla_u \Phi)^{-1}(0)\cap (\h \times [\lambda_-,\lambda_+])$ containing $(u_0,\lambda_+)$.
Since $C$ is bounded, there is $\xi\geq\gamma$ such that $C \subset B_{\xi}(\h)\times[\lambda_-,\lambda_+].$

Put in Lemma \ref{sep}
\begin{enumerate}
\item[(i)] $\ds K = \nabla_u \Phi^{-1}(0) \cap \left(cl(B_{\xi}(\h))\times[\lambda_-,\lambda_+]\right),$
\item[(ii)]$\ds A =C,$
\item[(iii)]$\ds B= \nabla_u \Phi^{-1}(0)  \cap \left(\partial
B_{\xi}(\h)\times[\lambda_-,\lambda_+]\right).$
\end{enumerate}

\nt Applying Lemma \ref{sep}, we obtain compact subsets $K_A,K_B \subset K$ such that $A\subset K_A,\; B\subset K_B,
\; K_A\cap K_B=\emptyset$ and $K_A\cup K_B= K$. Note that almost all $(u_{n_k},\lambda_+)\in K_B.$ Indeed,
$(u_{n_k},\lambda_+)\in C_{\gamma_{n_k}}$ and $C_{\gamma_{n_k}}\cap (\partial
B_{\gamma_{n_k}}(\h)\times[\lambda_-,\lambda_+])\neq \emptyset$. Hence $(u_{n_k},\lambda_+)\in K_B$ for all
$k\in\bN$ such that $\gamma_{n_k}\geq\xi$. Thus $(u_0,\lambda_+)$, as the limit of elements from the closed set
$K_B$, belongs to $K_B$. On the other hand, $(u_0,\lambda_+)\in A\subset K_A$ and $\dist(K_A,K_B)>0$, a
contradiction. We have just proved that $C \cap (\h \times [\lambda_-,\lambda_+])$ is unbounded.
\end{proof}

\br\label{rema}   Since $\nabla_u \Phi(\cdot,\lambda_{\pm})$ is of the form compact perturbation of the identity,
one can define a bifurcation index $\textsc{Bif}_{\mathrm{LS}}(\infty,[\lambda_-,\lambda_+]) \in \bZ$ as follows
$$\textsc{Bif}_{\mathrm{LS}}(\infty,[\lambda_-,\lambda_+])= \dls(\nabla_u \Phi(\cdot,\lambda_+),B_{\gamma}(\bH),0) - \dls(\nabla_u
\Phi(\cdot,\lambda_-),B_{\gamma}(\bH),0).$$ We realize that theorems similar to Theorem \ref{cont1} has been proved
for operators of the form compact perturbation of the identity (without gradient and equivariant structures), see
for instance Theorem 2.6 of \cite{[LESC]}.

However directly from the definition of the degree for $\sone$-equivariant gradient maps it follows that if
$\textsc{Bif}_{\mathrm{LS}}(\infty,[\lambda_-,\lambda_+]) \neq 0 \in \bZ$ then
$\textsc{Bif}(\infty,[\lambda_-,\lambda_+]) \neq \Theta \in U(\sone).$ On the other hand it can happen that
$\textsc{Bif}_{\mathrm{LS}}(\infty,[\lambda_-,\lambda_+]) = 0$ and $\textsc{Bif}(\infty,[\lambda_-,\lambda_+]) \neq
\Theta.$  \er

\bdf Let $C \subset \h \times \bR$ be  closed and connected. We say that  a symmetry breaking phenomenon for $C$
occurs if there are  $(u_0,\lambda_0) \in C$  and sequence $\{(u_n,\lambda_n)\} \subset C$ converging to
$(u_0,\lambda_0)$   such that $\sone_{u_n} \neq \sone_{u_0}$ for every $n \in \bN.$ \edf

\bco\label{wn1cont1} Let assumptions of Theorem \ref{cont1} be satisfied. Moreover, suppose that $ (\nabla_u
\Phi)^{-1}(0) \cap (\h^{\sone}\times [\lambda_-,\lambda_+])$ is bounded. Then,  there exists an unbounded closed
connected component $C$ of $(\nabla_u \Phi)^{-1}(0) \cap (\h \times [\lambda_-,\lambda_+])$ such that   the symmetry
breaking phenomenon for $C$ occurs or there exists at least one nontrivial solution of  equation \eqref{rbif}  such
that $(u,\lambda)\in (B_{\gamma}(\h)\times\{\lambda_-,\lambda_+\})\cap C$. \eco
\begin{proof}
By Theorem \ref{cont1} we obtain an unbounded  component $C$ of  $(\nabla_u \Phi)^{-1}(0) \cap (\h \times
[\lambda_-,\lambda_+])$ such that $C \cap \left(B_{\gamma}(\h)\times\{\lambda_-,\lambda_+\}\right)\neq\emptyset.$
Since $(\h^{\sone}\times [\lambda_-,\lambda_+])\cap (\nabla_u \Phi)^{-1}(0)$ is bounded, without loss of generality,
one can assume that
\begin{equation}\label{zeraogr}
(\nabla_u \Phi)^{-1}(0)\cap(\h^{\sone}\times [\lambda_-,\lambda_+])\subset B_{\gamma}(\h)\times
[\lambda_-,\lambda_+]
\end{equation}
Therefore the isotropy group of every element $u \in C \cap ((\h\setminus B_{\gamma}(\h))\times
[\lambda_-,\lambda_+])$ is different from $\sone.$ Thus, if $C \cap (\h^{\sone}\times
[\lambda_-,\lambda_+])\neq\emptyset,$ then  the symmetry breaking phenomenon for $C$  occurs. Otherwise
$C\subset\cN(\nabla_u \Phi)$ and  $C \cap \left(B_{\gamma}(\h)\times\{\lambda_-,\lambda_+\}\right)\neq\emptyset,$
which completes the proof.
\end{proof}

\br\label{uw1wn1cont1} Notice that if in Corollary \ref{wn1cont1} we have
$$(\nabla_u \Phi)^{-1}(0)  \cap (B_{\gamma}(\h)\times\{\lambda_{\pm}\})\subset\h^{\sone}\times\{\lambda_{\pm}\},$$ then
the symmetry breaking phenomenon for $C$ occurs. \er

\br\label{uw2wn1cont1} Notice that if in Corollary \ref{wn1cont1} we have $$(\nabla_u \Phi)^{-1}(0)\cap
(\h^{\sone}\times [\lambda_-,\lambda_+])=\{u_1,\dotsc,u_q \}\times [\lambda_-,\lambda_+]$$ and $\nabla_u^2
\Phi(u_i,\lambda)$ is an isomorphism for every $\lambda\in [\lambda_-,\lambda_+], i=1,\ldots,q,$ then $C
\subset\cN(\nabla_u \Phi)$. \er

\br\label{uw3wn1cont1} Under the assumptions of Corollary \ref{wn1cont1}. Since $$(\nabla_u
\Phi)^{-1}(0)\cap(\h^{\sone}\times [\lambda_-,\lambda_+])$$ is bounded, there is $\gamma > 0$ such that $$ (\nabla_u
\Phi)^{-1}(0)\cap(\h^{\sone}\times [\lambda_-,\lambda_+])\subset B_{\gamma}(\h)\times [\lambda_-,\lambda_+]. $$
Therefore we obtain
$$\dls((\nabla_u
\Phi(\cdot,\lambda_-))^{\sone},B_{\gamma}(\h)^{\sone},0)=\dls((\nabla_u
\Phi(\cdot,\lambda_+))^{\sone},B_{\gamma}(\h)^{\sone},0).$$ As a direct consequence of results due to Rabier
\cite{[RAB]} we obtain
\begin{equation}\label{rabier}
\dls(\nabla_u \Phi(\cdot,\lambda_{\pm}),B_{\gamma}(\h),0)=\dls((\nabla_u
\Phi(\cdot,\lambda_{\pm}))^{\sone},B_{\gamma}(\h)^{\sone},0).
\end{equation}
Summing up, we have obtained $\textsc{Bif}_{\mathrm{LS}}(\infty,[\lambda_-,\lambda_+])=0 \in \bZ.$ \er

\nt The following lemma is a parameterized extension of Corollary 3.1 of \cite{[FRR]}.

\bl\label{lem1cont1} Let $\Phi \in C^2_{\sone}(\h \times \bR,\bR)$ satisfy assumption \textrm{\bf(c1)}. Then for
every $(u_0,\lambda_0) \in (\nabla_u \Phi )^{-1}(0) \cap (\h^{\sone}\times \bR)$ there exist $\gamma > 0$ such that
if $(u,\lambda) \in (\nabla_u\Phi)^{-1}(0)\cap( B_{\gamma}(\h,u_0))\times(\lambda-\gamma,\lambda+\gamma))$, then
there exists $v \in \ker \nabla^2_u \Phi(u_0,\lambda_0)$ such that $\sone_u=\sone_v$. \el
\begin{proof}   Since  $\nabla^2_u \Phi(u_0,\lambda_0) :\h \rightarrow \h$ is a self-adjoint  Fredholm operator of
index $0$, we obtain $\h= \ker \nabla^2_u \Phi(u_0,\lambda_0) \oplus \im \nabla^2_u \Phi(u_0,\lambda_0)$. Let
$\pi:\h\rightarrow \ker \nabla^2_u \Phi(u_0,\lambda_0)$ and $Id - \pi:\h\rightarrow \im \nabla^2_u
\Phi(u_0,\lambda_0)$ stand for $\sone$-equivariant orthogonal projections.  Obviously
$$\nabla_u \Phi(u,\lambda)=0 \;\; \Leftrightarrow \;\;(\pi\circ\nabla_u \Phi)(u,\lambda)=0 \;\text{and}\; ((Id-\pi)\circ\nabla_u
\Phi)(u,\lambda)=0.$$ By the $\sone$-equivariant version of the implicit function theorem, we obtain that solutions
of $((Id-\pi)\circ\nabla_u \Phi)(u,\lambda)=0$ are of the form  $(v,\omega(v,\lambda),\lambda),$ where \linebreak
$v\in B_{\gamma}(\ker \nabla^2_u \Phi(u_0,\lambda_0),u_0),$ $ \lambda\in (\lambda_0-\gamma,\lambda_0+\gamma)$ for
sufficiently small $\gamma >0$ and $(v,\lambda) \rightarrow \omega(v,\lambda)$ is an $\sone$-equivariant
$C^1$-mapping.

Let $(u,\lambda) \in (\nabla_u\Phi)^{-1}(0)\cap( B_{\gamma}(\h,u_0))\times (\lambda_0-\gamma,\lambda_0+\gamma)$.
Therefore $(u,\lambda)=(v,\omega(v,\lambda),\lambda).$ Since  $\omega$ is $\sone$-equivariant, $\sone_{(v,\lambda)}
\subset \sone_{\omega(v,\lambda)}$ and consequently
$$\sone_u=\sone_{(u,\lambda)}=\sone_{(v,\omega(v,\lambda),\lambda)}=\sone_{(v,\lambda)}
\cap \sone_{\omega(v,\lambda)}=\sone_{(v,\lambda)}=\sone_{v}.$$
 \end{proof}

\nt As a direct consequence of Lemma \ref{lem1cont1} we obtain the following corollary.

\bco\label{wn2cont1} Let assumptions of Theorem \ref{cont1} be satisfied. Additionally, suppose that $\ker\nabla^2_u
\Phi(u,\lambda)\subset\h^{\sone}$ for every $u\in\h^{\sone}, \lambda\in [\lambda_-,\lambda_+].$  Then,
$$\textrm{either } C \subset \h^{\sone} \times [\lambda_-,\lambda_+]  \text{ or } C \subset \cN(\nabla_u \Phi).$$
If moreover $(\nabla_u \Phi)^{-1}(0)\cap B_{\gamma}(\h)\times\{\lambda_-,\lambda_+\}\subset \h^{\sone}\times
\{\lambda_-,\lambda_+\}$, then the symmetry breaking phenomenon for  $C$  does not occur. \eco
\begin{proof}
First of all notice that the set $C$ obtained by Theorem \ref{cont1} is closed and connected. Suppose, contrary to
our claim, that  $C \cap (\h^{\sone}\times [\lambda_-,\lambda_+])\neq\emptyset$ and $C\cap\cN(\nabla_u
\Phi)\neq\emptyset$. Then there exists $(u_0,\lambda_0) \in C \cap (\h^{\sone}\times [\lambda_-,\lambda_+])$ such
that in its any neighborhood  there exists an element $(u,\lambda)\in C \cap \cN(\nabla_u \Phi)$. Taking into
account that $\sone_u \neq \sone_{u_0}=\sone,$ the assumption and Lemma \ref{lem1cont1} we obtain a contradiction.
\end{proof}

\nt Let us put some additional assumptions on behaviour of the functional $\Phi$ at infinity. We would like to say
something more about behaviour of closed connected components of $(\nabla_u \Phi)^{-1}(0)$ at infinity. Suppose that
the functional $\Phi \in C^2_{\sone}(\h\times\bR,\bR)$ satisfies assumption \textrm{\bf(c1)} and the following
assumption:

\begin{enumerate}
\item[\textrm{\bf(c2)}] $\ds{\Phi(u,\lambda)=\frac{1}{2}\langle u,u\rangle_{\h}-\frac{1}{2} \langle K_{\infty}(\lambda)u ,u
\rangle_{\h} - \eta_{\infty}(u,\lambda),}$ where
\begin{enumerate}
\item[\bf (i)\rm] $K_{\infty}(\lambda): \h\rightarrow\h$ is a linear,
$\sone$-equivariant, self-adjoint, operator for every $\lambda\in\bR,$
\item[\bf (ii)\rm] the mapping $\h\times\bR\ni (u,\lambda) \mapsto
K_{\infty}(\lambda)u \in \h$ is compact,
\item[\bf (iii)\rm] $\nabla_u
\eta_{\infty}:\h\times\bR\rightarrow\h$ is a $\sone$-equivariant, compact operator such that $ \nabla_u
\eta_{\infty}(u,\lambda)=o(\| u\|),$\; as $\| u \|\rightarrow \infty$   uniformly on bounded $\lambda$-intervals.
\end{enumerate}
\end{enumerate}

\nt For $\lambda \in \bR$ define   $\nabla_u^2 \Phi(\infty,\lambda)= Id - K_{\infty}(\lambda).$ Fix arbitrary
$\lambda_0 \in \bR$ and assume that $\ker \nabla^2_u \Phi(\infty,\lambda_0) \neq \{0\}.$   Choose $\varepsilon>0,$
define $\lambda_{\pm}=\lambda_{0}\pm\varepsilon$ and assume that the following condition is fulfilled
\begin{equation}\label{izo}
\{\lambda\in [\lambda_-,\lambda_+] : \nabla^2_u \Phi(\infty,\lambda) \;\; \text{is not an
isomorphism}\}=\{\lambda_{0}\}.
\end{equation}

\nt It is easy to see that under the above assumptions  there exists $\gamma>0$ such that condition \eqref{zera} is
satisfied.

\bdf\label{meets} We say that an  unbounded closed connected set $C$ meets $(\infty,\lambda_{0})$, if for every
$\delta, \gamma > 0$
\begin{equation}\label{gghh}
C \cap \{(\h\setminus B_{\gamma}(\h)) \times [\lambda_{0}-\delta,\lambda_{0}+\delta]\} \neq \emptyset.
\end{equation}
\edf

\nt In the following theorem we localize points at which closed connected sets of solutions of equation \eqref{rbif}
meet infinity.

\bt\label{cont2} Let potential $\Phi \in C^2_{\sone}(\h\times\bR,\bR)$ satisfy assumption \textrm{\bf (c2)}. Choose
$\varepsilon, \gamma > 0, \; \lambda_{0},\lambda_{\pm}\in\bR$ such that \eqref{zera} and \eqref{izo} hold true. If
$\textsc{Bif}(\infty,[\lambda_-,\lambda_+])\neq\Theta\in U(\sone),$ then the statement of Theorem \ref{cont1} holds
true. Moreover, $C$ meets $(\infty,\lambda_{0})$. \et
\begin{proof}
The existence of an unbounded closed connected component $C$ of $\nabla_u \Phi^{-1}(0) \cap (\h \times
[\lambda_-,\lambda_+])$ satisfying $C \cap
\left(B_{\gamma}(\h)\times\{\lambda_{0}^-,\lambda_{0}^+\}\right)\neq\emptyset$, is a direct consequence of Theorem
\ref{cont1}. It remains to prove that $C$ meets $(\infty,\lambda_{0}).$ Note that it is sufficient to show, that
condition \eqref{gghh} holds true just for large $\varrho>0$ and small $\delta>0$. Choose any $\delta>0$ such that
$\delta<\varepsilon$. By assumption, for $\lambda\in [\lambda_{0}-\varepsilon,\lambda_{0}-\delta)\cup
(\lambda_{0}+\delta,\lambda_{0}+\varepsilon]$, $\nabla^2 \Phi(\infty,\lambda)$ is an isomorphism. Moreover, by
\textrm{\bf (c2)} we obtain  $\nabla_u \Phi(u,\lambda)=\nabla^2_u \Phi(\infty,\lambda)u +\nabla_u
\eta_{\infty}(u,\lambda),$ where $\nabla_u \eta_{\infty}(u,\lambda)=o(\| u\|),$\; as $\| u\|\rightarrow\infty$
uniformly on bounded $\lambda$-intervals i.e.
$$\forall_{\epsilon>0} \; \exists_{R_{\epsilon}>0} \;\forall_{\lambda\in [a,b]\subset\bR} \; \forall_{u\in\h}
\| u \| > R_{\epsilon} \Rightarrow
 \; \| \nabla_u
\eta_{\infty}(u,\lambda) \| < \epsilon \| u \|.$$

\nt Put $\ds{\epsilon=\frac{\| \nabla^2_u \Phi(\infty,\lambda)^{-1} \|^{-1}}{4}}$. Hence, for $\| u \|
>
 R_{\epsilon}$, we obtain
\begin{align*}
\| \nabla_u \Phi(u,\lambda) \| &= \| \nabla^2_u \Phi(\infty,\lambda)u +\nabla_u \eta_{\infty}(u,\lambda) \|
 \geq \| \nabla^2_u \Phi(\infty,\lambda)u \| -
\|\nabla_u \eta_{\infty}(u,\lambda) \|  \geq \\
& \geq  \frac{\| \nabla^2_u \Phi(\infty,\lambda)^{-1} \|^{-1}}{2} \| u \| - \frac{\| \nabla^2_u
\Phi(\infty,\lambda)^{-1}
\|^{-1}}{4} \| u \|  \geq \\
&\geq \frac{\| \nabla^2_u \Phi(\infty,\lambda)^{-1} \|^{-1}}{4} \| u \| >0.
\end{align*}

\nt Hence, for every $\varrho > R_{\epsilon}$,
$$C\cap((\h\setminus B_{\varrho}(\h,\infty))\times
[\lambda_{0}-\varepsilon,\lambda_{0}-\delta)\cup (\lambda_{0}+\delta,\lambda_{0}+\varepsilon])=\emptyset.$$ Since
$C$ is unbounded, $C \cap (\h\setminus B_{\varrho}(\h,\infty))\times [\lambda_{0}-\delta,\lambda_{0}+\delta] \neq
\emptyset,$ which completes the proof.
\end{proof}

\nt The principal significance of the lemma below is that it allows one to control the isotropy groups of solutions
of equation  \eqref{rbif}   sufficiently close to infinity.

\bl\label{lem1cont2} Let $\Phi \in C^2_{\sone}(\h \times \bR,\bR)$ satisfy assumption \textrm{\bf (c2)}. Then for
every $\lambda_0 \in \bR$ there exist $\gamma  > 0, \delta > 0$ such that if $(u,\lambda) \in
(\nabla_u\Phi)^{-1}(0)\cap(\h\setminus B_{\gamma}(\h))\times[\lambda_0-\delta,\lambda_0+\delta]$, then there exists
$v \in \ker(Id-K_{\infty}(\lambda_{0}))$ such that $\sone_u=\sone_v$.\el
\begin{proof}
Fix $\lambda_0 \in \bR$. By the $\sone$-equivariant version of the implicit function theorem at infinity (see
Theorem 3.2 of \cite{[FRR]}), we obtain that solutions of $\nabla_u \Phi(u,\lambda)=0$ in a neighborhood of
$(\infty,\lambda_0)$ are of the form $(v,\omega(v,\lambda),\lambda),$ where $ v\in \ker \nabla_u^2
\Phi(\infty,\lambda_0) \setminus cl(B_{\gamma}(\ker \nabla_u^2 \Phi(\infty,\lambda_0))),$ $ \lambda\in
[\lambda_{0}-\delta,\lambda_{0}+\delta]$ for some $\gamma,\delta
>0$ and the map
$(v,\lambda) \rightarrow \omega(v,\lambda)\in \im \nabla_u^2 \Phi(\infty,\lambda_0)$ is an $\sone$-equivariant
$C^1$-mapping. The rest of the proof is the same as the proof of Lemma \ref{lem1cont1}.
\end{proof}

\br\label{uw1lem1cont2} If moreover, assumptions of Theorem \ref{cont2} are satisfied, then without loss of
generality one can assume that $\delta \leq \varepsilon$. \er

\nt Below we  present some useful corollaries of Theorem \ref{cont2}. First of them is a counterpart of Corollary
\ref{cont1} at infinity, also based on Corollary 3.1 of \cite{[FRR]}.

\bco\label{wn2cont2} Let assumptions of Theorem \ref{cont2} be satisfied. Additionally suppose that $\ds{\ker\left(
(\nabla_u^2 \Phi(\infty,\lambda_{0}))\right) \cap \h^{\sone}=\{0\}}$. Then the statement of Theorem \ref{cont2}
holds true. Moreover, for closed connected set $C$ either phenomenon of symmetry breaking occurs or there exists at
least one nontrivial solution of  equation \eqref{rbif}  such that $(u,\lambda)\in C \cap (
B_{\gamma}(\h)\times\{\lambda_-,\lambda_+\})$. \eco
\begin{proof}
Note that by assumption and Lemma \ref{lem1cont2}, the isotropy group of any solution of  equation \eqref{rbif}
close to $(\infty,\lambda_{0})$ is different from $\sone$. Thus $(\h^{\sone}\times [\lambda_-,\lambda_+])\cap
(\nabla_u \Phi)^{-1}(0)$ is bounded and by Theorem \ref{cont2} and Corollary \ref{wn1cont1} the proof is completed.
\end{proof}

\bdf Let $\bV$ and $\bW$ be $\sone$-representations. We say that $\sone$-re\-pre\-sen\-ta\-tion $\bV$ is not
consistent with $\sone$-representation $\bW$, if $\sone_v \neq \sone_w$ for every $v \in \bV \setminus\{0\}, w \in
\bW \setminus\{0\}$. \edf

\br $\sone$-representation $\ds \bV = \bigoplus_{i=1}^{p} \bR[k_i,m_i]$ is not consistent with
$\sone$-representation $\ds \bW = \bigoplus_{j=1}^{q} \bR[k_j',m_j']$, if $\gcd(m_{i_1}',\dotsc,m_{i_r}')  \neq \gcd
( m_{j_1}',\dotsc,m_{j_s}'),$ for every \linebreak $\{i_1,\dotsc,i_r\}\subset\{1,\dotsc,p\},\;
\{j_1,\dotsc,j_s\}\subset\{1,\dotsc,q\}.$ \er


\bco \label{wn4cont2} Let assumptions of Theorem \ref{cont2} be satisfied. Additionally, suppose that
\begin{enumerate}
\item[(i)]  $(\nabla_u \Phi)^{-1}(0)\cap
(\h^{\sone}\times [\lambda_-,\lambda_+])=\{u_1,\dotsc,u_q \}\times
[\lambda_-,\lambda_+]$,
\item[(ii)]  $(\nabla_u \Phi)^{-1}(0)\cap
(B_{\gamma}(\h)\times \{\lambda_{\pm}\})=\{u_1,\dotsc,u_q \}\times \{\lambda_{\pm}\},$
\item[(iii)] $\{(u,\lambda)\in\{u_1,\dotsc,u_q \}\times
[\lambda_-,\lambda_+]\; : \;\nabla_u^2 \Phi(u,\lambda)\; \text{is not an isomorphism}\}=$ \linebreak
$\{(u_{i_1},\lambda_{i_1}),\dotsc,(u_{i_d},\lambda_{i_d})\}$,
\item[(iv)] $\ker(\nabla^2_u\Phi(u_{i_k},\lambda_{i_k}))$ is
not consistent with $\ker(\nabla^2_u\Phi(\infty,\lambda_{0}))$ for every $k=1,\dotsc,d.$
\end{enumerate}
Then the statement of Theorem \ref{cont2} holds true. Moreover, for $C$ phenomenon of symmetry breaking occurs. \eco
\begin{proof} By Theorem  \ref{cont2} we obtain an unbounded closed connected component $C$ of \linebreak
$(\nabla \Phi)^{-1}(0) \cap (\bH \times [\lambda_-,\lambda_+])$ such that $C \cap (B_{\gamma}(\bH) \times
\{\lambda_{\pm}\}) \neq \emptyset.$ From assumption (ii) it follows that $C \cap (\{u_1,\ldots,u_k\} \times
[\lambda_-,\lambda_+]) \neq \emptyset.$ Moreover, by assumption (iii) we obtain $C \cap (\{u_1,\ldots,u_k\} \times
[\lambda_-,\lambda_+]) \subset \{(u_{i_1},\lambda_{i_1}),\dotsc,(u_{i_d},\lambda_{i_d})\}.$ The rest of the proof is
a direct consequence of assumption (iv) and Lemmas \ref{lem1cont1}, \ref{lem1cont2}.
\end{proof}

\nt One can also proof the following slight generalization of Corollary \ref{wn4cont2}. Since the proof of the
following corollary is similar to the proof of Corollary \ref{wn4cont2} we omit it.

\bco \label{wn3cont2} Let assumptions of Theorem \ref{cont1} be satisfied. Additionally, suppose that
\begin{enumerate}
\item[(i)]  $\ds (\nabla_u \Phi)^{-1}(0)\cap
(\h^{\sone}\times [\lambda_-,\lambda_+])=\bigcup_{j=1}^q \{u_j\}\times [\lambda_-,\lambda_+]$,
\item[(ii)]  $\ds (\nabla_u \Phi)^{-1}(0)\cap
(B_{\gamma}(\h)\times \{\lambda_{\pm}\})=\bigcup_{j=1}^q \{u_j\}\times \{\lambda_{\pm}\},$
\item[(iii)] $\ds \{(u,\lambda)\in\{u_1,\dotsc,u_q \}\times
[\lambda_-,\lambda_+] : \ker \nabla_u^2  \Phi(u,\lambda) \neq \{0\}\}=\bigcup_{j=1}^d \{(u_{i_j},\lambda_{i_j})\},$
\item[(iv)] $\ds \{\lambda \in [\lambda_-,\lambda_+]\; : \ker \nabla_u^2 \Phi(\infty,\lambda) \neq \{0\}\}=
\bigcup_{j=1}^p \{\lambda_j^{\infty}\} \subset (\lambda_-,\lambda_+),$
\item[(v)] $\ker(\nabla^2_u\Phi(u_{i_k},\lambda_{i_k}))$ is
not consistent with $\ker(\nabla^2_u\Phi(\infty,\lambda_{j}^{\infty}))$ for every $k=1,\ldots,d$ and $j=1,\ldots,p.$
\end{enumerate}
Then the statement of Theorem \ref{cont2} holds true. Moreover,
\begin{enumerate}
  \item[a)] there is $j_0 \in \{1,\ldots,p\}$ such that $C$ meets $(\infty,\lambda_{j_0}),$
  \item[b)] for $C$ the phenomenon of symmetry breaking occurs.
\end{enumerate}
\eco


\numsec
\section{Connected Sets of   Periodic Solutions Bifurcating from Infinity}
\label{results}

In this section we study continuation of $2\pi$-periodic solutions of  family of autonomous second order Hamiltonian systems of the form
\begin{equation}\label{ns}
\big(E_{\lambda}\big) \:\:\:\: \;\;\;   \begin{cases}
    \ddot{u}(t) =  -  \nabla_u V(u(t),\lambda),&  \\
    u(0)=u(2\pi), & \\
    \dot{u}(0)=\dot{u}(2\pi),
  \end{cases}
\end{equation}
where
\begin{enumerate}
  \item[\textbf{(a1)}] $V \in C^2(\bR^n \times \bR,\bR),$
  \item[\textbf{(a2)}] $\ds V(x,\lambda)=\frac{1}{2} (A(\lambda)x,x)+\eta(x,\lambda),$ where $(\cdot, \cdot)$ is the usual
  scalar product in  $\bR^n.$
  \item[\textbf{(a3)}] $A(\lambda)$ is real symmetric matrix for every $\lambda\in\bR,$
  \item[\textbf{(a4)}] $\nabla_x \eta(x,\lambda)=o(\|x\|),$ as $ \|x\| \rightarrow\infty$ uniformly  on bounded $\lambda$-intervals.
\end{enumerate}

Define a separable Hilbert space
\\\centerline{$\h^1_{2\pi} = \{u : [0,2\pi] \rightarrow \bR^n : \text{ u is abs. cont., } u(0)=u(2\pi), \dot u \in
L^2([0,2\pi],\bR^n)\}$} \\
with a scalar product given by the formula $\ds \langle u,v\rangle_{\h^1_{2\pi}} = \int_0^{2\pi} (\dot u(t), \dot
v(t)) + (u(t),v(t)) \; dt.$ The space $\left(\h^1_{2\pi} ,\langle\cdot,\cdot\rangle_{\h^1_{2\pi}}\right)$ is an
orthogonal $\sone$-representation with the $\sone$-action given by shift in time.

It is well known that solutions of  system \eqref{ns}  are in one to one correspondence with critical points of an $\sone$-invariant
$C^2$-functional $\Phi_V: \h^1_{2\pi} \times \bR \rightarrow \bR$ given by the formula
\begin{equation}\label{nsf}
\Phi_V(u,\lambda) = \frac{1}{2} \int_0^{2\pi} \mid \dot u(t) \mid^2 \; dt - \int_0^{2\pi} V(u(t),\lambda) \; dt.
\end{equation}
\nt Moreover, it is known that $\nabla_u^2 \Phi_V(\infty,\lambda)= Id - L_{A(\lambda)}$, where $L_{A(\lambda)}:\h^1_{2\pi} \rightarrow\h^1_{2\pi}$
is a linear, self-adjoint, $\sone$-equivariant and compact operator defined by the formula $\langle L_{A(\lambda)} (u),v\rangle_{\h^1_{2\pi}}=\ds
\int_0^{2\pi} (u(t)+A(\lambda)u(t),v(t))dt.$ By Corollary 5.1.1. of \cite{[FRR]}, $\nabla_u^2 \Phi_V(\infty,\lambda)$ is an isomorphism iff $\ds
\sigma(A(\lambda)) \cap \left\{k^2 : k \in \bN \cup \{0\}\right\}= \emptyset.$ Note that $\Phi_V:\h^1_{2\pi} \times \bR \rightarrow \bR$ satisfies
assumptions \textrm{\bf (c1), (c2)} of the previous section.

\nt Let us put two additional assumptions:

\nt \textbf{(a5)} assume that there exist $\lambda_-,\lambda_+>0$ such that the set of solutions of $(E_{\lambda_{\pm}})$ is bounded in
$\bH^1_{2\pi}$, i.e. there exists $\gamma>0$ such that

\begin{equation}\label{izolinfty}
(\nabla_u\Phi_V(\cdot,\lambda_{\pm})^{-1}(0)\cap((\h^1_{2\pi} \setminus
B_{\gamma}(\h^1_{2\pi}))\times\{\lambda_{\pm}\})=\emptyset,
\end{equation}

\nt \textbf{(a6)} assume that
\begin{itemize}
\item $\ds \sigma(A(\lambda_-)) \cap \left\{ k^2
: \ k \in \bN \right\} = \left\{ (k_1^-)^2, \ldots,(k_r^-)^2\right\},$
\item $\ds \sigma(A(\lambda_+)) \cap \left\{ k^2
: \ k \in \bN \right\} = \left\{(k_1^{+})^2, \ldots,(k_s^{+})^2\right\}.$
\end{itemize}

\nt Put
$$\ds \bK =\bigcup_{\{i_1, \dotsc, i_l\} \in \{1, \dotsc, r\}} \{\gcd(k^-_{i_1},\dotsc,k^-_{i_l})\}
 \cup \bigcup_{\{i_1, \dotsc, i_m\} \in \{1, \dotsc, s\}}
\{\gcd(k^+_{i_1},\dotsc,k^+_{i_m})\}.$$

\nt If $\sigma(A(\lambda_{\pm})) \cap \left\{ k^2 : \ k \in \bN \right\} =\emptyset$, then it is understood that
$\bK=\emptyset.$ For $\alpha \in \bR$  we will denote by $\mu_A(\alpha)$ the multiplicity of $\alpha$ considered as
an eigenvalue of matrix $A.$ If $\alpha \notin \sigma(A) $ then it is understood that $\mu_A(\alpha)=0.$ For every
$k \in \bN \cup \{0\}$ define
\begin{enumerate}
  \item $\ds \sigma_k(A,2\pi)=\sigma(A) \cap\left(k^2,+\infty \right),$
  \item $\ds j_k(A,2\pi)=\sum_{\alpha \in \sigma_k(A,2\pi)} \mu_A(\alpha).$
\end{enumerate}
Put $\ds \ind(-\nabla_x V(\cdot,\lambda_{\pm}),\infty)=\lim_{\alpha \rightarrow \infty}\deg_{\rm B}(-\nabla_x
V(\cdot,\lambda_{\pm}),B_{\alpha}(\bR^n,0),0),$ where $\deg_{\rm B}$ denotes the Brouwer degree.

\bt\label{eqcont1} Let assumptions \textrm{\bf(a1)}-\textrm{\bf(a6)} be satisfied. Additionally, suppose that one of the following conditions
holds:
\begin{enumerate}
\item[(i)] $\mathrm{ind}(\nabla_x V(\cdot,\lambda_+),\infty)\neq\mathrm{ind}(
\nabla_x V(\cdot,\lambda_-),\infty),$
\item[(ii)] $\mathrm{ind}(\nabla_x V(\cdot,\lambda_+),\infty)=\mathrm{ind}(
\nabla_x V(\cdot,\lambda_-),\infty)\neq 0$ and there exists $k\in\bN\setminus\bK$ such that $j_k\left(A(\lambda_+),2\pi\right)\neq
j_k\left(A(\lambda_-),2\pi\right).$
\end{enumerate}
Then there exists an unbounded closed connected component $C \subset \bH^1_{2\pi}\times[\lambda_-,\lambda_+]$ of
solutions of  system \eqref{ns} such that $C
\cap\left(B_{\gamma}(\h^1_{2\pi})\times\{\lambda_-,\lambda_+\}\right)\neq\emptyset$. \et
\begin{proof}
First of all notice  that $\Phi_V: \h^1_{2\pi} \times \bR \rightarrow \bR$ given by formula \eqref{nsf} satisfies condition \textrm{\bf (c1)}.
\begin{enumerate}
\item[(i)] By Lemma 5.2.3. of \cite{[FRR]}, $$\dg_{\sone}(\nabla_u
\Phi_V(\cdot,\lambda_{\pm}),B_{\gamma}(\h^1_{2\pi}))=\mathrm{ind}(-\nabla_x V(\cdot,\lambda_{\pm}),\infty).$$ That
is why we obtain
$$\textsc{Bif}_{\sone}(\infty,[\lambda_-,\lambda_+]) =$$
$$=\dg_{\sone}(\nabla_u \Phi_V(\cdot,\lambda_+),B_{\gamma}(\h^1_{2\pi}))
-\dg_{\sone}(\nabla_u \Phi_V(\cdot,\lambda_-),B_{\gamma}(\h^1_{2\pi}))=$$
$$=\mathrm{ind}(-\nabla_x V(\cdot,\lambda_+),\infty) -\mathrm{ind}(-\nabla_x V(\cdot,\lambda_-),\infty)\neq 0.$$
\item[(ii)] By Lemma 5.2.3. of \cite{[FRR]},
$$\dg_{\bZ_k}(\nabla_u
\Phi_V(\cdot,\lambda_{\pm}),B_{\gamma}(\h^1_{2\pi}))=\mathrm{ind}(-\nabla_x V(\cdot,\lambda_{\pm}),\infty) \cdot
j_k\left(A(\lambda_{\pm}),2\pi\right).$$ Therefore we have
$$\textsc{Bif}_{\bZ_k}(\infty,[\lambda_-,\lambda_+]) =$$
$$=\dg_{\bZ_k}(\nabla_u \Phi_V(\cdot,\lambda_+),B_{\gamma}(\h^1_{2\pi}))
-\dg_{\bZ_k}(\nabla_u \Phi_V(\cdot,\lambda_-),B_{\gamma}(\h^1_{2\pi}))=$$
$$=\mathrm{ind}(-\nabla_x V(\cdot,\lambda_+),\infty) \cdot
j_k\left(A(\lambda_+),2\pi\right)-\mathrm{ind}(- \nabla_x V(\cdot,\lambda_-),\infty) \cdot j_k\left(A(\lambda_-),2\pi\right)\neq 0.$$
\end{enumerate}
Since $\textsc{Bif}(\infty,[\lambda_-,\lambda_+])\neq\Theta \in U(\sone),$ the rest of the proof is a direct
consequence of Theorem \ref{cont1}.
\end{proof}

\bt \label{faktzeraogr}   Theorem \ref{eqcont1} remains true if the  assumption \textrm{\bf(a5)} is replaced by
\begin{enumerate}
  \item[(a)] $\bK=\emptyset,$
  \item[(b)] $(\nabla_x V(\cdot,\lambda_{\pm})^{-1}(0)\cap((\bR^n\setminus B_{\gamma}(\bR^n))\times\{\lambda_{\pm}\})=\emptyset.$
\end{enumerate}
\et
\begin{proof} Notice that $(\h^1_{2\pi})^{\sone}=\bR[n,0]$ and that
$\nabla_u\Phi_V(\cdot,\lambda_{\pm})^{\sone}=-  \nabla_x V(\cdot,\lambda_{\pm})$. From  Lemma \ref{lem1cont2} it
follows that  for every $(u,\lambda)\in(\nabla_u\Phi_V)^{-1}(0)$ close to $(\infty,\lambda_{\pm})$, there exists
$v\in\ker\nabla_u^2 \Phi_V(\infty,\lambda_{\pm})$ such that $\sone_{(u,\lambda)}=\sone_v$. Combining  the
assumptions with Lemma 5.1.1 and Corollary 5.1.1. of \cite{[FRR]} we obtain that $\ker\nabla_u^2
\Phi_V(\infty,\lambda_{\pm})\subset\bR[n,0].$ Therefore $\sone_{(u,\lambda)}=\sone$ and
$\nabla_u\Phi_V(u,\lambda)=0$ iff $\nabla_u V(u,\lambda)=0$.
\end{proof}

\bdf We say that $2\pi \geq T>0$ is a period of function $u\in\h^1_{2\pi}$ if $u(t+T)=u(t)$ for every $t\in
[0,2\pi].$ We say that $T_{\min}\geq 0$ is a minimal period of function $u\in\h^1_{2\pi}$ if $T_{\min}=\inf \{T>0:
u(t+T)=u(t) \text{ for every } t \in [0,2\pi]\}$. \edf

\br Notice that if $u\in(\h^1_{2\pi})^{\sone}$, i.e. $u=\text{const}$, $T_{\min}=0$ and therefore $T_{\min}$ is not a period of function $u$.
Nevertheless, we call $T_{\min}=0$ the minimal period of a constant function $u.$ \er

\bco\label{wn1eqcont1} Let assumptions of Theorem \ref{eqcont1} be satisfied. If additionally $(\nabla_x
V)^{-1}(0)\cap (\bR^n\times[\lambda_-,\lambda_+])$ is bounded, then conclusion of   Theorem \ref{eqcont1} holds
true. Moreover, continuum $C$ emanates from the set of stationary solutions and contains solutions with different
minimal periods or there exists at least one non-stationary solution $(u,\lambda)$ of  system \eqref{ns} such that
$(u,\lambda)\in (B_{\gamma}(\h^1_{2\pi})\times\{\lambda_-,\lambda_+\})\cap C$.\eco

\begin{proof}
Note that $(\h^1_{2\pi})^{\sone}=\bR[n,0].$ It is clear that solutions with different isotropy group have different minimal periods. Since all the
assumptions of Corollary \ref{wn1cont1} are satisfied, we obtain our assertion.
\end{proof}

\br Under assumptions of Corollary \ref{wn1eqcont1}, if moreover equations $(E_{\lambda_{\pm}})$ possesses only
stationary periodic solutions then continuum $C$ contains solutions with different minimal periods. \er

\bco\label{wn2eqcont1} Let assumptions of Theorem \ref{eqcont1} be satisfied. Additionally, suppose that $\ker
\nabla_u^2\Phi_V(u,\lambda)\subset (\bH^1_{2\pi})^{\sone}=\bR[n,0]$ for every $u  \in (\bH^1_{2\pi})^{\sone}$ and
$\lambda\in[\lambda_-,\lambda_+],$ then conclusion of   Theorem \ref{eqcont1} holds true. Moreover, either $C
\subset  (\bH^1_{2\pi})^{\sone} \times [\lambda_-,\lambda_+]$ or $C$ contains only non-stationary solutions. If
additionally  equations $(E_{\lambda_{\pm}})$ possesses only stationary periodic solutions   then $C$ consists of
stationary solutions of  system \eqref{ns} . \eco
\begin{proof}
Immediate consequence of Corollary \ref{wn2cont1}.
\end{proof}

\nt Let us put the following assumption

\nt \textrm{\bf (a7)} fix $\lambda_0\in\bR$ and choose  $\lambda_- < \lambda_+$    such that
\begin{equation}\label{izom}
\left\{\lambda\in [\lambda_-,\lambda_+] : \sigma(A(\lambda)) \cap \left\{k^2 : k \in \bN \cup \{0\}\right\}\neq \emptyset\right\}=\left\{\lambda_0
\right\}.
\end{equation}

\nt Combining assumption \eqref{izom} with Corollary 5.1.1 of \cite{[FRR]} we obtain that $\nabla_u^2 \Phi(\infty,\lambda_{\pm}) : \bH^1_{2\pi}
\rightarrow \bH^1_{2\pi}$ is a linear isomorphism. Therefore assumption \textrm{\bf (a5)} is satisfied.

\bt\label{eqcont2} Let assumptions \textrm{\bf(a1)}-\textrm{\bf(a4)}, \textrm{\bf(a7)}   be satisfied. Additionally, suppose that at least one of
the following conditions holds:
\begin{enumerate}
\item[(i)]
$(-1)^{j_0(A(\lambda_+),2\pi)}\neq(-1)^{j_0(A(\lambda_-),2\pi)}$,
\item[(ii)] there exists $k\in \bN$ such that
$j_k\left(A(\lambda_+),2\pi\right)\neq j_k\left(A(\lambda_-),2\pi\right)$.
\end{enumerate}
Then there exists an unbounded closed connected component $C \subset \bH^1_{2\pi} \times [\lambda_-,\lambda_+]$ of
solutions of  system \eqref{ns} such that $C \cap\left(B_{\gamma}(\h^1_{2\pi})\times\{\lambda_-,\lambda_+\}\right)
\neq \emptyset$. Moreover, $C$ meets $(\infty,\lambda_0)$. \et

\begin{proof}
Note that $\Phi_V: \h^1_{2\pi} \times \bR \rightarrow \bR$ given by formula \eqref{nsf} satisfies \textrm{\bf (c2)}.
\begin{enumerate}
\item[(i)] By Lemma 5.2.2. and Remark 5.2.2. of \cite{[FRR]},
$$\dg_{\sone}(\nabla_u
\Phi_V(\cdot,\lambda_{\pm}),B_{\gamma}(\h^1_{2\pi}))=(-1)^{j_0(A(\lambda_{\pm}),2\pi)}.$$ Therefore we obtain
$$\textsc{Bif}_{\sone}(\infty,[\lambda_-,\lambda_+]) =$$
$$=\dg_{\sone}(\nabla_u \Phi_V(\cdot,\lambda_+),B_{\gamma}(\h^1_{2\pi}))
-\dg_{\sone}(\nabla_u \Phi_V(\cdot,\lambda_-),B_{\gamma}(\h^1_{2\pi}))=$$
$$=(-1)^{j_0(A(\lambda_+),2\pi)} -(-1)^{j_0(A(\lambda_-),2\pi)}\neq 0.$$
\item[(ii)] By Lemma 5.2.2. and Remark 5.2.2. of \cite{[FRR]},
$$\dg_{\bZ_k}(\nabla_u
\Phi_V(\cdot,\lambda_{\pm}),B_{\gamma}(\h^1_{2\pi}))=(-1)^{j_0(A(\lambda_{\pm}),2\pi)} \cdot
j_k\left(A(\lambda_{\pm}),2\pi\right).$$ That is why we have
$$\textsc{Bif}_{\bZ_k}(\infty,[\lambda_-,\lambda_+]) =$$
$$=\dg_{\bZ_k}(\nabla_u \Phi_V(\cdot,\lambda_+),B_{\gamma}(\h^1_{2\pi}))
-\dg_{\bZ_k}(\nabla_u \Phi_V(\cdot,\lambda_-),B_{\gamma}(\h^1_{2\pi}))=$$
$$=(-1)^{j_0(A(\lambda_+),2\pi)} \cdot
j_k\left(A(\lambda_+),2\pi\right)-(-1)^{j_0(A(\lambda_-),2\pi)} \cdot j_k\left(A(\lambda_-),2\pi\right).$$
\end{enumerate}
Summing up, $\textsc{Bif}(\infty,[\lambda_-,\lambda_+])\neq\Theta$. The rest of the proof is a direct consequence of
Theorem \ref{cont2}.
\end{proof}

\nt Recall that by Corollary 5.1.2. of \cite{[FRR]} $$\ker \nabla_u^2 \Phi_V(\infty,\lambda_0)= \ker(Id -
L_{A(\lambda_0)}) \approx \bigoplus_{k=0}^{\infty} \bR\left[\mu_{A(\lambda_0)}\left( k^2\right),k\right].$$  Note
that for almost every  $k\in\bN\cup\{0\}$, $k^2\notin\sigma(A(\lambda_0))$ and hence
$\mu_{A(\lambda_0)}\left(k^2\right)=0$. Since $\bR[0,k]=\{0\}$, $\dim \ker \nabla_u^2 \Phi_V(\infty,\lambda_0) <
\infty.$

\bco\label{wn1eqcont2} Let assumptions of Theorem \ref{eqcont2} be satisfied. Suppose that
$$\ds \sigma(A(\lambda_0))\cap\left\{k^2 : k\in\bN\cup\{0\}\right\}=
\left\{k_0^2,k_1^2,\ldots,k_r^2\right\},$$ where $0 \leq k_0 < k_1 \dotsc < k_r.$
\begin{enumerate}
\item[(i)] If $\det A(\lambda_0)=0$, then for every solution $(u,\lambda)$ of  system \eqref{ns}  in $\bH^1_{2\pi} \times
[\lambda_-,\lambda_+]$ sufficiently close to $(\infty,\lambda_0)$ its minimal period $T_{\min}$ is equal to zero (u=const) or to
$\frac{2\pi}{gcd(k_{i_1}, \dotsc, k_{i_s})}$ for some $\{k_{i_1},\dotsc,k_{i_s}\}\subset\{k_1,\dotsc,k_r\}$.
\item[(ii)] If $\det A(\lambda_0)\neq 0$, then for every solution $(u,\lambda)$ of  system \eqref{ns}  in $\bH^1_{2\pi} \times
[\lambda_-,\lambda_+]$ sufficiently close to $(\infty,\lambda_0)$ its minimal period $T_{\min}$ is equal to $\frac{2\pi}{gcd(k_{i_1}, \dotsc,
k_{i_s})}$ for some $\{k_{i_1},\dotsc,k_{i_s}\}\subset\{k_0,\dotsc,k_r\}$.
\end{enumerate}
\eco
\begin{proof} By assumption and Corollary 5.1.2. of \cite{[FRR]} we have
$$\ker \nabla_u^2 \Phi_V(\infty,\lambda_0)= \ker(Id -
L_{A(\lambda_0)})  \approx \bigoplus_{i=0}^{r} \bR\left[\mu_{A(\lambda_0)}\left(k_i^2\right),k_i\right].$$
By Lemma \ref{lem1cont2} any solution
$(u,\lambda)$ of  system \eqref{ns} sufficiently close  to $(\infty,\lambda_0)$ has the same isotropy group as  some element of $\ker \nabla_u^2
\Phi_V(\infty,\lambda_0).$ Therefore if $\det A(\lambda_0)=0$, then the possible isotropy group of any solution is equal to $\sone$ or
$\bZ_{gcd(k_{i_1}, \dotsc, k_{i_s})}$ for some $\{k_{i_1},\dotsc,k_{i_s}\}\subset\{k_1,\dotsc,k_r\},$ which completes the proof of (i). Otherwise,
it is equal to $\bZ_{gcd(k_{i_1}, \dotsc, k_{i_s})}$ for some $\{k_{i_1},\dotsc,k_{i_s}\}\subset\{k_0,\dotsc,k_r\}$, which completes the proof of
(ii).
\end{proof}

\bco\label{wn2eqcont2} Let assumptions of Theorem \ref{eqcont2} be satisfied. If additionally $\det A(\lambda_0)\neq
0$ then conclusion of Theorem \ref{eqcont2} holds true. Moreover, continuum $C$ emanates from the set of stationary
solutions and contains solutions with different minimal periods or there exists at least one non-stationary solution
such that $(u,\lambda)\in (B_{\gamma}(\h^1_{2\pi})\times\{\lambda_-,\lambda_+\})\cap C$. \eco
\begin{proof}
By Lemma 5.1.1 and Corollary 5.1.1 we obtain $\ker \nabla_u^2 \Phi_V(\infty,\lambda_0) \cap (\bH^1_{2\pi})^{\sone}=\{0\}$ iff $\det
A(\lambda_0)\neq 0$. The rest of the proof is a direct consequence of Lemma \ref{lem1cont2}.\end{proof}

From now on we consider  special case of   system \eqref{ns}. Namely, we consider system
\begin{equation}\label{nss}
   \begin{cases}
    \ddot{u}(t) =  -  \lambda^2 \nabla V(u(t)),&  \\
    u(0)=u(2\pi), & \\
    \dot{u}(0)=\dot{u}(2\pi),
  \end{cases}
\end{equation}
where
\begin{enumerate}
\item[\textbf{(b1)}] $V \in C^2(\bR^n,\bR),$
\item[\textbf{(b2)}] $V(x)=\frac 12 (A x,x)+\eta(x),$
\item[\textbf{(b3)}] $A$ is a real symmetric matrix,
\item[\textbf{(b4)}] $\nabla \eta(x)=o(\|x\|),$  as $\|x\| \rightarrow \infty,$
\item[\textbf{(b5)}] $(\nabla V)^{-1}(0)$ is bounded,
\item[\textbf{(b6)}] $\mathrm{ind}(\nabla V,\infty)\neq 0.$
\end{enumerate}

It is easy to show that $\nabla^2_u \Phi_V(\infty,\lambda)$ is not an isomorphism if and only if \bc$\lambda \in
\left\{\frac{k}{\sqrt{\alpha}}: k\in\bN,\; \alpha\in\sigma_{+}(A)\right\} \text{ or } \det A \neq 0.$\ec

\bl \label{lemind}  Fix $k_0 \in \bN, \alpha_0 \in \sigma_{+}(A)$ and choose $\lambda_- < \lambda_+$ such that \bc $\left[\lambda_-,\lambda_+
\right]\cap\left\{\frac{k}{\sqrt{\alpha}}:\; k\in\bN,\;\alpha\in\sigma_{+}(A) \right\}=\left\{\frac{k_0}{\sqrt{\alpha_0}}\right\}.$\ec Then
$\textsc{Bif}(\infty,[\lambda_-,\lambda_+]) \in U(\sone)$ is well-defined. Moreover,
$$\textsc{Bif}_{\bZ_{k_0}}(\infty,[\lambda_-,\lambda_+])=\mathrm{ind}(-\nabla V,\infty)\cdot\mu_A(\alpha_0).$$\el
\begin{proof} Since $(\nabla_u \Phi(\cdot,\lambda_{\pm}))^{-1}(0) \subset \bH^1_{2\pi}$ is bounded, $\textsc{Bif}(\infty,[\lambda_-,\lambda_+]) \in
U(\sone)$ is well-defined. Applying Lemma 5.2.2 of \cite{[FRR]}, we obtain:
$$\textsc{Bif}_{\bZ_{k_0}}(\infty,[\lambda_-,\lambda_+])=$$
$$=\dg_{\bZ_{k_0}}(Id-L_{\lambda_+^2A},B_{\gamma}(\h^1_{2\pi}))
-\dg_{\bZ_{k_0}}(Id-L_{\lambda_-^2 A},B_{\gamma}(\h^1_{2\pi}))=$$
$$=\mathrm{ind}(-\lambda_+^2 \nabla V,\infty) \cdot j_{k_0}\left(\lambda_+^2 A,2\pi\right)-\mathrm{ind}(-\lambda_-^2
\nabla V,\infty) \cdot j_{k_0}\left(\lambda_-^2 A,2\pi\right)=$$
$$=\mathrm{ind}(-\nabla V,\infty) \cdot
\left(j_{k_0}\left(\lambda_+^2 A,2\pi\right)-j_{k_0}\left(\lambda_-^2 A,2\pi\right)\right)=$$
$$=\mathrm{ind}(-\nabla V,\infty) \cdot \left(\sum_{\alpha \in \sigma_{k_0}(\lambda_+^2 A,2\pi)} \mu_{(\lambda_+^2
A)}(\alpha)-\sum_{\alpha \in \sigma_{k_0}(\lambda_-^2 A,2\pi)} \mu_{(\lambda_-^2 A)}(\alpha)\right)=$$
$$=\mathrm{ind}(-\nabla V,\infty) \cdot \mu_A(\alpha_0).$$
\end{proof}

\nt The following theorem is a consequence of Theorem \ref{eqcont2}.

\bt\label{eqcont3} Let assumptions \textrm{\bf (b1)}-\textrm{\bf (b6)} be fulfilled. Then for every $$\lambda_0\in
\left\{\frac{k}{\sqrt{\alpha}}: k\in\bN,\; \alpha\in\sigma_{+}(A)\right\}$$ there exists an unbounded closed
connected component $C(\lambda_0) \subset \bH^1_{2\pi} \times [\lambda_-,\lambda_+]$ of solutions of system
\eqref{nss}   such that $C(\lambda_0) \cap
\left(B_{\gamma}(\h^1_{2\pi})\times\{\lambda_-,\lambda_+\}\right)\neq\emptyset,$ where $\lambda_- < \lambda_+$
satisfy $\left[\lambda_-,\lambda_+ \right]\cap\left\{\frac{k}{\sqrt{\alpha}}:\; k\in\bN,\;\alpha\in\sigma_{+}(A)
\right\}=\left\{\lambda_0\right\}.$ Moreover, $C(\lambda_0)$ meets $(\infty,\lambda_0)$.\et

\nt Fix $\lambda_0=\frac{k_0}{\sqrt{\alpha_0}}$ for some $k_0\in\bN,\; \alpha_0\in \sigma_{+}(A)$.

\bco\label{wn1eqcont3} Let assumptions of Theorem \ref{eqcont3} be satisfied. Assume additionally that
\begin{enumerate}
\item[(i)] $(\nabla V)^{-1}(0)=\{u_1,\dotsc,u_q\}$,
\item[(ii)] the only periodic solutions of $\big(E_{\lambda_{\pm}^2}\big)$ are the critical points of $V,$
\item[(iii)] $\{(u,\lambda)\in\{u_1,\dotsc,u_q \}\times
[\lambda_-,\lambda_+]\; : \;\ds \sigma(\lambda^2 \nabla^2 V(u_i)) \cap \left\{k^2 : k \in \bN \cup \{0\}\right\} \neq \emptyset\}=$ \linebreak
$\{(u_{i_1},\lambda_{i_1}),\dotsc,(u_{i_d},\lambda_{i_d})\}$,
\item[(iv)] $\ker(\nabla^2_u\Phi_V(u_{i_k},\lambda_{i_k}))$ is
not consistent with $\ker(\nabla^2_u\Phi_V(\infty,\lambda_0))$ for all $k=1,\dotsc,d.$
\end{enumerate}
Then there exists an unbounded closed connected component $C(\lambda_0) \subset \bH^1_{2\pi} \times
[\lambda_-,\lambda_+]$ of solutions of system \eqref{nss}  such that $C(\lambda_0) \cap
\left(B_{\gamma}(\h^1_{2\pi})\times\{\lambda_-,\lambda_+\}\right)\neq\emptyset$ and  $C(\lambda_0)$ meets
$(\infty,\lambda_0)$. Moreover, $C(\lambda_0)$ contains solutions with different minimal periods.  \eco
\begin{proof} Note that $\ds \sigma(\lambda^2 \nabla^2 V(u_{i_k})) \cap \left\{k^2
: k \in \bN \cup \{0\}\right\}= \emptyset$ implies that$\nabla_u^2 \Phi_V(u_{i_k},\lambda)$ is an isomorphism for
every $k=1,\dotsc,d.$ Therefore applying Corollary \ref{wn4cont2} we complete the proof.
\end{proof}




\numsec
\section{Examples}
\label{examples}

In this section we discuss three examples of potentials in order to illustrate   results proved in the previous
section. We consider   system \eqref{ns}  with simple potential $V$ and show that assumptions of our theorems are
satisfied.

\bex  Define potential  $ V: \bR^n\times\bR \rightarrow \bR$ as follows
\begin{equation}\label{spp}
V (x,\lambda)= \frac{1}{2} (A(\lambda) x, x) + W (x,\lambda) =\frac{1}{2} (A(\lambda) x, x) +
\frac{-\lambda^2}{\sqrt{\|x\|^2 +a}},
\end{equation}
where $a > 0$ and  $A(\lambda)$ is a real symmetric $(n\times n)$-matrix for every $\lambda\in\bR$. Consider system
\eqref{ns}  with potential \eqref{spp}. Put $n=4, a=1,$ $ \lambda_{\pm}=\pm 1$ and define
$$A(\lambda)=\left[\begin{array}{cccc}\ds
 \lambda^2-1 & 0 & 0 & 0 \\
  0 & \sqrt{2}+\lambda & 0 & 0 \\
  0 & 0 & \lambda-\sqrt{2} & 0 \\
  0 & 0 & 0 & \sqrt{5}+\lambda
\end{array} \right].$$
Systems $\big(E_{\pm 1}\big)$ are resonant at infinity because
\begin{equation}\label{www}
\sigma(A(\pm 1)) \cap \left\{ k^2 : k \in \bN \cup \{0\}\right\} = \{0\}.
\end{equation}

\nt Notice that assumptions \textrm{\bf (a1)}-\textrm{\bf (a4), (a6)} of Theorem   \ref{eqcont1} are satisfied. Moreover,
\begin{enumerate}
\item $(\nabla_x V(\cdot,\pm 1))^{-1}(0)$ is bounded because $\#(\nabla_x V(\cdot,\pm 1))^{-1}(0) < \infty$
  (consequence of Lem\-ma 6.2 of \cite{[FRR]}),
\item $\bK=\emptyset$    (consequence of \eqref{www}).
\end{enumerate}
Applying Theorem \ref{faktzeraogr} we show that assumption \textrm{\bf (a5)} of Theorem \ref{eqcont1} is fulfilled.

Moreover,
\begin{enumerate}
  \item $\ind(-\nabla_x V(\cdot,\pm 1),\infty)=(-1)^{n-m^-(A(\pm 1))}=(-1)^{4-1}=-1$ (consequence of Lemma 6.4 of
  \cite{[FRR]}),
  \item $j_1(A(+1),2\pi)=2 \neq 1 = j_1(A(-1),2\pi).$
\end{enumerate}
Applying Theorem \ref{eqcont1} we obtain an unbounded closed connected component $C \subset \bH^1_{2\pi} \times
[-1,+1]$ of solutions of  system \eqref{ns}  such that $C \cap\left(B_{\gamma}(\h^1_{2\pi})\times\{-1,+1\}\right)
\neq \emptyset.$

Additionally, taking into consideration that
\begin{enumerate}
\item $(\nabla_x V(\cdot,\pm 1))^{-1}(0)$ is bounded,
\item  $\{\lambda \in (-1,+1) : \sigma(A(\lambda)) \cap \left\{ k^2 : k \in \bN \cup
\{0\}\right\} \neq \emptyset\}=\{\lambda_0=1-\sqrt{2}\},$
\item $\sigma(A(1-\sqrt{2}))=\{1\},$
\end{enumerate}
and Corollary \ref{wn1eqcont2} we obtain that the continuum meets $(\infty,1-\sqrt{2})$ and that any solution
$(u,\lambda) \in C$ of  system \eqref{ns}  sufficiently close to $(\infty,1-\sqrt{2})$ has minimal period equal to
$2\pi.$ \eex

\bex Define potential  $ V: \bR^n\times\bR \rightarrow \bR$ as follows
\begin{equation}\label{sitnikovpotential}
V (x,\lambda)= \frac{1}{2} (A(\lambda) x, x) + W (x,\lambda) =\frac{1}{2} (A(\lambda) x, x) +
\frac{-1}{\sqrt{\|x\|^2 +a}},
\end{equation}
where $a > 0$ and  $A(\lambda)$ is a real symmetric $(n\times n)$-matrix for every $\lambda\in\bR$.

\nt Consider   system \eqref{ns}  with potential \eqref{sitnikovpotential}. Put $n=4, a=1$ and define
$$A(\lambda)=\left[\begin{array}{cccc}\ds
 4 + \lambda& 0 & 0 & 0 \\
  0 & 2 & 0 & 0 \\
  0 & 0 & 2 & 0 \\
  0 & 0 & 0 & 2
\end{array} \right].$$
System $\big(E_{0}\big)$ is resonant at infinity because
\begin{equation}\label{rer}
\sigma(A(0)) \cap \left\{ k^2 : k \in \bN \cup \{0\}\right\} = \{4\}.
\end{equation}
Moreover, put $\lambda_{\pm}=\pm (1\slash 2)$ and notice that
\begin{equation}\label{rer1}
\sigma(A(\lambda)) \cap \left\{ k^2 : k \in \bN \cup \{0\}\right\}=\emptyset
\end{equation}
for every $\lambda \in [-1 \slash 2,+1 \slash 2] \setminus \{0\}.$

\nt Since $j_2(A(\frac{1}{2}))=1 \neq 0= j_2(A(-\frac{1}{2})),$ all the assumptions of Theorem \ref{eqcont2} are
fulfilled. Therefore there exists an unbounded closed connected component $C \subset \bH^1_{2\pi} \times [-1 \slash
2,+1 \slash 2]$ of solutions of system \eqref{ns}  such that $C \cap\left(B_{\gamma}(\h^1_{2\pi})\times\{-1 \slash
2,+1 \slash 2\}\right) \neq \emptyset$ and that $C$ meets $(\infty,0)$.

\nt Properties of potential $V$ have been precisely studied in \cite{[FRR]}. \nt Stationary solutions of  system
\eqref{ns}  have the following properties:
\begin{enumerate}
  \item $(\nabla_x V)^{-1}(0) \cap (\bR^4 \times [-1 \slash 2,+1 \slash 2])=\{0\} \times [-1 \slash 2,+1 \slash 2]$
  (consequence of Lemma 6.2 of \cite{[FRR]}),
  \item $\nabla^2_{xx}(0,\lambda)=A(\lambda)+Id,$ for every $\lambda \in [-1 \slash 2,+1 \slash 2]$
  (consequence of Lemma 6.1 of \cite{[FRR]}),
  \item $\sigma(\nabla^2_{xx} V(0,\lambda)) \cap \left\{ k^2 : k \in \bN \cup
\{0\}\right\}=\emptyset$ for every $\lambda \in [-1 \slash 2,+1 \slash 2]$ (consequence of (2)).
\end{enumerate}

\nt Moreover, by \eqref{rer}, \eqref{rer1}  and Corollary   \ref{wn1eqcont2} we obtain that  any solution
$(u,\lambda) \in C$ of  system \eqref{ns}  sufficiently close to $(\infty,0)$ has minimal period equal to $\pi.$
Additionally, from {\rm (3)} and Remark \ref{uw2wn1cont1} it follows that continuum $C$ consist of non-stationary
solutions.\eex

\bex  Consider   system \eqref{ns}  with potential \eqref{sitnikovpotential}.  Put $n=5, a=1, \lambda_{\pm}=\pm 1$ and define
$$A(\lambda)=\left[\begin{array}{ccccc}\ds
 4+\frac{\lambda^2}{2} & 0 & 0 & 0 & 0 \\
  0 & \lambda^3-\sqrt{10} & 0 & 0 \\
  0 & 0 & \ds 9+\frac{\lambda^2}{2} & 0 & 0 \\
  0 & 0 & 0 & \lambda^3+\sqrt{10} & 0 \\
  0 & 0 & 0 & 0 & \ds 25+\frac{\lambda^2}{2}
\end{array} \right].$$

\nt It is easy to see that
\begin{enumerate}
\item $\sigma(A(\lambda)) \cap \left\{ k^2 : k \in \bN \cup \{0\}\right\}= \emptyset$ for every
$\lambda\in[-1,1] \setminus \{0\},$
\item $\sigma(A(0)) \cap \left\{ k^2 : k \in \bN \cup \{0\}\right\}= \{4,9,25\}.$
\end{enumerate}
Hence assumptions  \textrm{\bf (a1)}-\textrm{\bf (a4), (a7)} of Theorem \ref{eqcont2} are fulfilled.

\nt Since $j_2(A(1),2\pi)=4 \neq 3= j_2(A(-1),2\pi),$ all the assumption of Theorem \ref{eqcont2} are satisfied.
Therefore there exists an unbounded closed connected component $C$ of solutions of system \eqref{ns} in
$\bH^1_{2\pi}\times[-1,1]$ such that $C \cap\left(B_{\gamma}(\h^1_{2\pi})\times\{-1,1\}\right)\neq\emptyset$ and
$C$ meets $(\infty,0)$. Moreover, by $\mathrm{(2)}$ and Corollary \ref{wn1eqcont2} (ii) any solution $(u,\lambda)
\in C$ sufficiently close to $(\infty,0)$ possesses the minimal period $T_{min}\in\left\{\ds 2\pi, \ds \pi,
\ds\frac{2\pi}{3}, \ds\frac{2\pi}{5}\right\}.$ \eex

\end{document}